\documentclass[12pt]{article}
\usepackage{amssymb}
\usepackage{amsmath}
\usepackage{tikz}
\newtheorem{Theorem}{Theorem}[section]

\newtheorem{Lemma}[Theorem]{Lemma}

\newtheorem{Corollary}[Theorem]{Corollary}
\newcommand{\qed}{\hfill $\Box$}
\newenvironment{Proof}{ \emph{Proof.}}{\qed}
\def\bc{\begin{center}}
\def\ec{\end{center}}

\begin{document}
\title{\bf On the variety generated by all semirings of order two}
\author{Aifa Wang, Lili Wang\\
   {\small School of Mathematical Sciences, Chongqing University of Technology}\\
 {\small Chongqing, 400054, P.R. China}\\
{\small wangaf@cqut.edu.cn}\\
}
\date{}
\maketitle \vskip -4pt \baselineskip 16pt \noindent
\textbf{ABSTRACT.} 
There are ten distinct two-element semirings up to isomorphism, denoted \( L_2, R_2, M_2, D_2, N_2, T_2, Z_2, W_2, Z_7 \), and \( Z_8 \) (see \cite{bk}). Among these, the multiplicative reductions of \( M_2, D_2, W_2 \), and \( Z_8 \) form semilattices, while the additive reductions of \( L_2, R_2, M_2, D_2, N_2 \), and \( T_2 \) are idempotent semilattices, commonly referred to as \emph{idempotent semirings}. 
In 2015, Vechtomov and Petrov \cite{vp} studied the variety generated by \( M_2, D_2, W_2 \), and \( Z_8 \), proving that it is finitely based. In the same year, Shao and Ren \cite{srii} examined the variety generated by the six idempotent semirings, demonstrating that every subvariety of this variety is finitely based.

This paper systematically investigates the variety generated by all ten two-element semirings. We prove that this variety contains exactly 480 subvarieties, each of which is finitely based.



\section{ Introduction}

A \emph{variety} is a class of algebras of the same type that is closed under subalgebras, homomorphic images, and direct products. Birkhoff's theorem asserts that a class of algebras is a variety if and only if it is an equational class. One of the central questions in the study of varieties is the finite basis problem, which asks whether a variety can be defined by a finite set of identities. If this is the case, the variety is said to be \emph{finitely based}; otherwise, it is \emph{nonfinitely based}. An algebra \( A \) is termed finitely based (or nonfinitely based) if the variety generated by \( A \) is finitely based (or nonfinitely based).

The concept of a finite basis has led to significant research in algebraic theory. In 1951, Lyndon \cite{lci} demonstrated that all two-element algebras are finitely based and posed the problem of whether every finite algebra is finitely based. This question was later answered negatively by Lyndon \cite{lcii}, who showed that a certain 7-element groupoid is nonfinitely based.

Over the years, various finite classical algebras have been shown to be finitely based. For example, every finite group \cite{oats}, every finite associative ring \cite{krrl, liv}, every finite lattice \cite{mcnr}, and every commutative semigroup \cite{perp} are finitely based. However, not all finite semigroups and semirings are finitely based. The first example of a nonfinitely based finite semigroup was provided by Perkins \cite{perp}, and Dolinka \cite{dlk} gave the first example of a nonfinitely based finite semiring.

To address the finite basis problem for finite algebras, Tarski \cite{tasa} posed the question: Is there an algorithm to decide whether a finite algebra is finitely based? McKenzie \cite{mker} negatively answered this question for finite groupoids. However, the problem remains open when restricted to finite semigroups and semirings.

Inspired by these issues, several authors have studied the finite basis problem for semirings of finite (small) order. For instance, Guzmán \cite{gf} investigated the variety generated by the two-element distributive lattice and the two-element finite field, proving that it is finitely based. Ghosh, Pastijn, and Zhao \cite{gho} studied the variety generated by all \( \mathbf{ai} \)-semirings of order two, whose multiplicative reduct is a band, and demonstrated that this variety is finitely based. Shao and Ren \cite{srii} considered the variety \( \mathbf{aiSr}(2) \) generated by all \( \mathbf{ai} \)-semirings of order two, proving that it is finitely based. Similarly, Vechtomov and Petrov \cite{vp} studied the variety generated by all semirings of order two whose multiplicative reduct is a semilattice, showing that this variety is finitely based.

On the other hand, Zhao, Ren, Crvenković, Shao, and Dapić \cite{rzcsd} proved that, except for the semiring \( S_7 \), all three-element \( \mathbf{ai} \)-semirings are finitely based. Recently, Jackson, Ren, and Zhao \cite{jrz} showed that the semiring \( S_7 \) is nonfinitely based. McNulty and Willard \cite{mgfw} initiated the study of the finite basis problem for algebras of order three. The present paper continues this line of investigation by systematically exploring the finite basis problem for the semiring variety generated by all ten semirings of order two.

The article is organized as follows:
\begin{itemize}
  \item \textbf{Section 2}: Introduces the basic definitions and notation related to semirings, which will form the foundation for our main results.
  \item \textbf{Section 3}: Discusses the variety \( \mathbf{Sr}(2) \) generated by all semirings of order two. We prove that \( \mathbf{Sr}(2) \) is finitely based and provide a finite basis for it.
  \item \textbf{Section 4}: Presents our main results, where we establish a mapping from \( \mathbf{Sr}(2) \) to \( \mathbf{aiSr}(2) \), effectively dividing the subvariety lattice into 64 intervals. Using the characterization of these intervals (Theorem \ref{3.5}), we derive the subvariety lattice of \( \mathbf{Sr}(2) \) (Theorem \ref{3.6}).
\end{itemize}

\section{Preliminaries}

By a \emph{semiring}, we mean an algebra \( (S, +, \cdot) \) such that:

\begin{itemize}
    \item \( (S, +) \) is a commutative semigroup;
    \item \( (S, \cdot) \) is a semigroup;
    \item The distributive laws \( x(y + z) \approx xy + xz \) and \( (y + z)x \approx yx + zx \) hold in \( S \).
\end{itemize}

A semiring is a generalization of both rings and distributive lattices. Semirings are widely used in algebra and computer science, as evidenced by works such as \cite{glk,gol}. It is well-known that the class of all semirings forms a variety, denoted by \( \mathbf{Sr} \).

Let \( \mathbf{V} \) be a variety, \( \mathcal{L}(\mathbf{V}) \) be the lattice of subvarieties of \( \mathbf{V} \), and \( \textrm{Id}_{\mathbf{V}}(X) \) denote the set of all identities satisfied by \( \mathbf{V} \). A variety \( \mathbf{V} \) is said to be \emph{hereditarily finitely based} if all members of \( \mathcal{L}(\mathbf{V}) \) are finitely based. In particular, if a variety \( \mathbf{V} \) is finitely based and \( \mathcal{L}(\mathbf{V}) \) forms a finite lattice, then \( \mathbf{V} \) is hereditarily finitely based \cite{rz}.

Up to isomorphism, there are exactly ten semirings of order two, denoted by \( L_2, R_2, M_2,\) \( D_2, N_2, T_2, Z_2, W_2, Z_7, \) and \( Z_8 \), respectively (see \cite{bk}). Each of these semirings has the underlying set \( \{0, 1\} \). The addition and multiplication operation tables for these semirings are given in the following table 1.

\begin{table}[h!]
\centering
\renewcommand{\arraystretch}{1.2} % Increase row height
\begin{tabular}{l|c|c|l}
\hline
\textbf{Semiring} & \textbf{Addition} & \textbf{Multiplication} & \textbf{Equational Basis} \\
\hline
$L_2$ & 
\begin{tabular}{cc}
   0 & 1 \\
   1 & 1 \\
\end{tabular} & 
\begin{tabular}{cc}
   0 & 0 \\
   1 & 1 \\
\end{tabular} & 
\( x + x \approx x \), \( xy \approx x \) \\
\hline
$R_2$ & 
\begin{tabular}{cc}
   0 & 1 \\
   1 & 1 \\
\end{tabular} & 
\begin{tabular}{cc}
   0 & 1 \\
   0 & 1 \\
\end{tabular} & 
\( x + x \approx x \), \( xy \approx y \) \\
\hline
$M_2$ & 
\begin{tabular}{cc}
   0 & 1 \\
   1 & 1 \\
\end{tabular} & 
\begin{tabular}{cc}
  0 & 1 \\
  1 & 1 \\
\end{tabular} & 
\( x + x \approx x \), \( x + y \approx xy \) \\
\hline
$D_2$ & 
\begin{tabular}{cc}
   0 & 1 \\
   1 & 1 \\
\end{tabular} & 
\begin{tabular}{cc}
   0 & 0 \\
   0 & 1 \\
\end{tabular} & 
\( x + x \approx x \), \( x^2 \approx x \), \( xy \approx yx \), \( x + xy \approx x \) \\
\hline
$N_2$ & 
\begin{tabular}{cc}
   0 & 1 \\
   1 & 1 \\
\end{tabular} & 
\begin{tabular}{cc}
   0 & 0 \\
   0 & 0 \\
\end{tabular} & 
\( x + x \approx x \), \( xy \approx zt \), \( z + xy \approx z \) \\
\hline
$T_2$ & 
\begin{tabular}{cc}
   0 & 1 \\
   1 & 1 \\
\end{tabular} & 
\begin{tabular}{cc}
   1 & 1 \\
   1 & 1 \\
\end{tabular} & 
\( x + x \approx x \), \( xy \approx zt \), \( z + xy \approx xy \) \\
\hline
$Z_2$ & 
\begin{tabular}{cc}
   0 & 0 \\
   0 & 0 \\
\end{tabular} & 
\begin{tabular}{cc}
   0 & 0 \\
   0 & 0 \\
\end{tabular} & 
\( x + y \approx z + u \), \( xy \approx x + y \) \\
\hline
$W_2$ & 
\begin{tabular}{cc}
   0 & 0 \\
   0 & 0 \\
\end{tabular} & 
\begin{tabular}{cc}
   0 & 0 \\
   0 & 1 \\
\end{tabular} & 
\( x + y \approx z + u \), \( x^2 \approx x \), \( xy \approx yx \) \\
\hline
$Z_7$ & 
\begin{tabular}{cc}
   0 & 1 \\
   1 & 0 \\
\end{tabular} & 
\begin{tabular}{cc}
   0 & 0 \\
   0 & 0 \\
\end{tabular} & 
\( x + x + y \approx y \), \( xy \approx x + x \) \\
\hline
$Z_8$ & 
\begin{tabular}{cc}
   0 & 1 \\
   1 & 0 \\
\end{tabular} & 
\begin{tabular}{cc}
   0 & 0 \\
   0 & 1 \\
\end{tabular} & 
\( x + x + y \approx y \), \( x^2 \approx x \), \( xy \approx yx \) \\
\hline
\end{tabular}
\caption{All two-element semirings}
\end{table}

The equational basis for each of these semirings, which generates the corresponding semiring variety, has been provided in previous works. For instance, the bases for \( L_2, R_2, M_2, D_2, N_2, \) and \( T_2 \) can be found in \cite{gho}, while the bases for \( Z_2, W_2, Z_7, \) and \( Z_8 \) are derived from \cite{259}.

Let \( X \) be a fixed countably infinite set of variables, and \( X^{+} \) the free semigroup on \( X \). We now introduce some notational conventions.

Let \( u \approx v \) be a semiring identity, where \( u = u_1 + \cdots + u_m \) and \( v = v_1 + \cdots + v_\ell \), with \( u_i, v_j \in X^{+} \). Let \( q \in X^{+} \). We define the following notations:

\begin{itemize}
  \item The \emph{content} of \( q \), denoted \( c(q) \), is the set of variables occurring in \( q \).
  \item The \emph{head} of \( q \), denoted \( h(q) \), is the first variable occurring in \( q \).
  \item The \emph{tail} of \( q \), denoted \( t(q) \), is the last variable occurring in \( q \).
  \item The \emph{length} of \( q \), denoted \( |q| \), is the number of variables occurring in \( q \), where each variable is counted as many times as it occurs in \( q \).
  \item The \emph{content} of \( u \), denoted \( C(u) \), is the set \( \bigcup \{ c(u_i) \mid i = 1, 2, \dots, m \} \).
  \item The \emph{head} of \( u \), denoted \( H(u) \), is the set \( \bigcup \{ h(u_i) \mid i = 1, 2, \dots, m \} \).
  \item The \emph{tail} of \( u \), denoted \( T(u) \), is the set \( \bigcup \{ t(u_i) \mid i = 1, 2, \dots, m \} \).
  \item The number of times \( u_i \) occurs in \( u \), denoted \( \text{occ}(u_i, u) \).
\end{itemize}

Using the above definitions, we can state the following lemma:

\begin{Lemma}\label{1.1}
Let \( u \approx v \) be any nontrivial semiring identity, where \( u = u_1 + \cdots + u_m \) and \( v = v_1 + \cdots + v_\ell \), with \( u_i, v_j \in X^+ \), and \( 1 \leq i \leq m, 1 \leq j \leq \ell \). Then:

\begin{enumerate}
  \item[$(1)$] \( L_2 \models u \approx v \) if and only if \( H(u) = H(v) \);
  \item[$(2)$] \( R_2 \models u \approx v \) if and only if \( T(u) = T(v) \);
  \item[$(3)$] \( M_2 \models u \approx v \) if and only if \( C(u) = C(v) \);
  \item[$(4)$] \( D_2 \models u \approx v \) if and only if for all \( u_i \in u, v_k \in v \), there exist \( u_l \in u \) and \( v_j \in v \) such that \( c(v_j) \subseteq c(u_i) \) and \( c(u_l) \subseteq c(v_k) \);
  \item[$(5)$] \( N_2 \models u \approx v \) if and only if \( \{ u_i \in u \mid |u_i| = 1 \} = \{ v_i \in v \mid |v_i| = 1 \} \);
  \item[$(6)$] \( T_2 \models u \approx v \) if and only if \( \{ u_i \in u \mid |u_i| \geq 2 \} \neq \emptyset \) and \( \{ v_i \in v \mid |v_i| \geq 2 \} \neq \emptyset \);
  \item[$(7)$] \( Z_2 \models u \approx v \) if and only if \( \forall x \in X, \, u \neq x \) and \( v \neq x \);
  \item[$(8)$] \( W_2 \models u \approx v \) if and only if \( m = n = 1 \) and \( c(u_1) = c(v_1) \), or \( m, n \geq 2 \);
  \item[$(9)$] \( Z_7 \models u \approx v \) if and only if \( \{ u_i \mid |u_i| = 1, \, \text{occ}(u_i, u) \text{ is odd} \} = \{ v_j \mid |v_j| = 1, \, \text{occ}(v_j, v) \text{ is odd} \} \);
  \item[$(10)$] \( Z_8 \models u \approx v \) if and only if for every non-empty \( A \subseteq X \), the number of \( u_i \in u \) such that \( c(u_i) \subseteq A \) is odd if and only if the number of \( v_j \in v \) such that \( c(v_j) \subseteq A \) is odd.
\end{enumerate}
\end{Lemma}

The results in Lemma \ref{1.1} follow directly from the equational bases of the semirings \( L_2, R_2, M_2, D_2, N_2, T_2, Z_2, W_2, Z_7, \) and \( Z_8 \) in Table 1.

By Lemma \ref{1.1}, we can state the following lemma, with the proof omitted:

\begin{Lemma}\label{1.2}
The variety \( \mathbf{Sr}(2) \) satisfies the following identities:

\begin{eqnarray}
x^2 y & \approx & xy \label{721}, \\
xy^2 & \approx & xy \label{271}, \\
xyzt & \approx & xzyt \label{2}, \\
(xy)^2 & \approx & xy \label{3}, \\
x + yz & \approx & x + yz + 2xz \label{223}, \\
x + yz & \approx & x + yz + 2yx \label{233}, \\
x + yz & \approx & x + yz + 2x^2 \label{723}, \\
x + yz & \approx & x + yz + 2xyz \label{244}, \\
x + yz & \approx & x + yz + 2yzx \label{252}, \\
x + y & \approx & x + 3y \label{162}.
\end{eqnarray}
\end{Lemma}

For further notations and terminologies used in this paper, the reader is referred to \cite{bur}, \cite{gho}, \cite{hw}, and \cite{jmh}.

\section{Equational basis of  $\mathbf{Sr}(2)$}

\begin{Theorem}\label{31312}
The semiring variety $\mathbf{Sr}(2)$ is determined by the identities (\ref{721})--(\ref{162}).
\end{Theorem}

\begin{Proof}
By Lemma \ref{1.2}, $\mathbf{Sr}(2)$ satisfies the identities (\ref{721})--(\ref{162}). To prove that every identity $u \approx v$ satisfied by $\mathbf{Sr}(2)$, where $u = u_1 + \cdots + u_m$ and $v = v_1 + \cdots + v_n$ with $u_i, v_j \in X^+$, can be derived from the identities satisfied by $\mathbf{Sr}$ and the identities (\ref{721})--(\ref{162}), we consider the following cases.

\textbf{Case 1.} \( m = n = 1 \) and \( C(u) = C(v) \). We consider two subcases:

\begin{itemize}
    \item If \( u_1 = x \), then since \( N_2 \models u \approx v \), it follows that \( v_1 = x \), so \( u \equiv v = x \).
    
    \item If \( u_1 = x_1^{i_1} x_2^{i_2} \cdots x_m^{i_m} \) with \( m \geq 2 \), then, since \( L_2, R_2, M_2, T_2 \models u \approx v \), we have
    \[
    v_1 = y_1^{i_1} y_2^{i_2} \cdots y_m^{i_m},
    \]
    where 
    \[
    \{x_1, x_2, \dots, x_m\} = \{y_1, y_2, \dots, y_m\}.
    \]
    Applying identities (\ref{721}), (\ref{271}), and (\ref{2}), we derive \( u_1 \approx v_1 \).
\end{itemize}

\textbf{Case 2.} \( m, n \geq 2 \). We distinguish three subcases:

\begin{itemize}
    \item If \( \ell(u_i) = 1 \) and \( \ell(v_j) = 1 \) for all \( u_i \in u \), \( v_j \in v \), then from \( M_2, N_2, Z_8 \models u \approx v \) and identity (\ref{162}), it follows that \( u \approx v \).
    
    \item If \( \ell(u_i) \geq 2 \) and \( \ell(v_j) \geq 2 \) for all \( u_i \in u \), \( v_j \in v \), then for some \( v_j = x_1 \cdots x_\ell \), \( \ell \geq 2 \), note that \( \mathbf{Sr}(2) \models u \approx u + 2v_j \). By \( L_2 \models u \approx u + 2v_j \), there exists a monomial in \( u \) of the form \( x_1 t_1^{(1)} \), with \( t_1^{(1)} \in X^* \). Similarly, \( u \) contains a monomial \( t_2^{(1)} x_\ell \), with \( t_2^{(1)} \in X^* \). Dually, there exist monomials of the form:
    \[
    t_3^{(1)} x_2 t_4^{(1)}, \quad \dots, \quad t_{2\ell-3}^{(1)} x_{\ell-1} t_{2\ell-2}^{(1)}, \quad t_{2\ell-1}^{(1)} \in X^*.
    \]
    Moreover, by \( M_2 \models u \approx u + 2v_j \), there exists a monomial \( u_s \in u \) such that \( c(u_s) \subseteq c(v_j) \); and by \( D_2 \models u \approx u + 2v_j \), there exists a monomial \( u_f \in u \) with \( |u_f| \geq 2 \).

    Thus, we can rewrite \( u \) as:
    \[
    u = u' + u_s + u_f + x_1 t_1^{(1)} + t_3^{(1)} x_2 t_4^{(1)} + \cdots + t_{2\ell-3}^{(1)} x_{\ell-1} t_{2\ell-2}^{(1)} + t_2^{(1)} x_\ell.
    \]
    Applying identity (\ref{723}), we obtain:
    \[
    u \approx u' + u_f + u_s + 2u_s^2 + x_1 t_1^{(1)} + \cdots + t_2^{(1)} x_\ell.
    \]
    Using identities (\ref{223})--(\ref{252}), we further derive:
    \[
    u \approx u' + u_f + u_s + 2u_s^2 + x_1 t_1^{(1)} + \cdots + t_2^{(1)} x_\ell + 2x_1 u_s^2 + \cdots + 2x_\ell u_s^2.
    \]
    Applying identities (\ref{244}) and (\ref{252}) yields:
    \[
    u \approx u' + u_f + u_s + 2u_s^2 + x_1 t_1^{(1)} + \cdots + t_2^{(1)} x_\ell + 2x_1 u_s^2 + \cdots + 2x_1 u_s^4 x_2 u_s^2 + \cdots.
    \]
    Finally, using identity (\ref{162}), we conclude:
    \[
    u \approx u + 2v_j.
    \]
    Similarly, we obtain:
    \[
    v \approx v + 2u_i.
    \]
    Therefore,
    \[
    u \approx u + 2v \approx v + 2u \approx v.
    \]

    \item If there exist both \( u_k \) with \( \ell(u_k) = 1 \) and \( u_\ell \) with \( \ell(u_\ell) \geq 2 \), and likewise for \( v \), write
    \[
    u = u^{(1)} + u^{(2)}, \quad \text{where } u^{(1)} = u_1' + \cdots + u_k',\; |u_i'| = 1;\quad u^{(2)} = u_1'' + \cdots + u_\ell'',\; |u_j''| \geq 2.
    \]
    Similarly, write
    \[
    v = v^{(1)} + v^{(2)}, \quad \text{where } v^{(1)} = v_1' + \cdots + v_r',\; |v_p'| = 1;\quad v^{(2)} = v_1'' + \cdots + v_h'',\; |v_q''| \geq 2.
    \]
    Using similar arguments as before, we derive:
    \[
    u^{(1)} \approx v^{(1)}, \quad u^{(2)} \approx u^{(2)} + 2v_j^{(2)}, \quad v^{(2)} \approx v^{(2)} + 2u_i^{(2)}.
    \]
    Hence,
    \[
    u \approx u^{(1)} + u^{(2)} \approx v^{(1)} + u^{(2)} + 2v^{(2)} \approx v^{(1)} + v^{(2)} + 2u^{(2)} \approx v^{(1)} + v^{(2)} \approx v.
    \]
\end{itemize}
\end{Proof}

\begin{Corollary}\label{3131}The semiring variety $\mathbf{Sr}(2)$ is determined by the identities (\ref{2})--(\ref{233}), (\ref{162}).
\end{Corollary}
\begin{Proof} First, note that since
\[
x^2x\overset{(\ref{3})}\approx(x^2x)^2=x^6\overset{(\ref{3})}\approx x^2,
\]
the identity (4) implies
\begin{eqnarray}
x^3\approx x^2,\label{14}
\end{eqnarray}
Then, using identities (3), (4) and (11), we have:
\[xy\overset{(\ref{3})}\approx (xy)^2\overset{(\ref{2})}\approx x^2y^2 \overset{(\ref{14})}\approx x^3y^2 \overset{(\ref{2})}\approx x(xy)^2\overset{(\ref{3})}\approx x(xy).
\]
This shows that the identities (3), (4), (11) imply the identity (1). Similarly, the identities (3), (4), (11) imply the identity (2).

On the other hand, we compute the following using identities (\ref{223}), (\ref{233}), and (\ref{721}):
\[
x + yz \overset{(\ref{223})} \approx x + yz + 2xz \overset{(\ref{233})} \approx yx + x + 2xz + 2x^2 \overset{(\ref{223})} \approx x + yz + 2x^2,
\]
\[
x + yz \overset{(\ref{721})} \approx x + yyz \overset{(\ref{223})} \approx x + yyz + 2xyz \overset{(\ref{721})} \approx x + yz + 2xyz,
\]
\[
x + yz \overset{(\ref{271})} \approx x + yzz \overset{(\ref{233})} \approx x + yzz + 2yzx \overset{(\ref{271})} \approx x + yz + 2yzx.
\]
Thus, we conclude that \( \mathbf{Sr}(2) \) is determined by the identities (\ref{2})--(\ref{233}), (\ref{162}).
\end{Proof}

In 2015, Shao and Ren \cite{srii} studied the semiring variety $\mathbf{aiSr}(2)$ generated by the six idempotent semirings of order two and obtained the following result.
\begin{Lemma}\cite{srii}\label{2.1}
 $\mathbf{aiSr}(2)$ is determined by the identities (\ref{2}), (\ref{3}), (\ref{1}), and (\ref{4}), and ${\cal L}(\mathbf{aiSr}(2))$ is a 64-element Boolean lattice.
\begin{eqnarray}
2x & \approx & x; \label{1}\\
x+yz & \approx & x+yz+xz+yx.\label{4}
\end{eqnarray}
\end{Lemma}
From Lemma \ref{31312} and \ref{2.1}, the following Corollary is immediate.
\begin{Corollary}\label{2.3}$\mathbf{aiSr}(2)$ is the subvariety of $\mathbf{Sr}(2)$ determined by the identity (\ref{1}).
\end{Corollary}

\section{ The lattice of subvarieties of  $\mathbf{Sr}(2)$}

Throughout this section, let \( t(x_1,\ldots,x_n) \) denote a term that contains no variables other than \( x_1,\ldots,x_n \) (though not necessarily all of them). Let \( S \in \mathbf{Sr}(2) \), and define
\[
E^+(S) = \{a \in S \mid 2a = a\},
\]
where each element of \( E^+(S) \) is called an \emph{additive idempotent} of \( (S,+) \). Note that \( \mathbf{Sr}(2) \) satisfies the following identities:
\begin{align}
(2x)+(2y) &\approx 2(x+y); \label{19} \\
(2x)(2y) &\approx 2xy. \label{19'}
\end{align}

Moreover, by identities (\ref{19}) and (\ref{19'}), it follows that \( E^+(S) = \{2a \mid a \in S\} \) forms a subsemiring of \( S \).

To describe the lattice \( {\cal L}(\mathbf{Sr}(2)) \), consider the mapping
\begin{equation}\label{20}
\varphi: {\cal L}(\mathbf{Sr}(2)) \rightarrow {\cal L}(\mathbf{aiSr}(2)),\quad \mathbf{V} \mapsto \mathbf{V} \cap \mathbf{aiSr}(2).
\end{equation}
It is easy to verify that for each \( \mathbf{V} \in {\cal L}(\mathbf{Sr}(2)) \), we have
\[
\varphi(\mathbf{V}) = \{E^+(S) \mid S \in \mathbf{V}\}.
\]

Let \( \mathbf{V} \subseteq \mathbf{aiSr}(2) \) be the subvariety defined by the identities
\[
u^{(\ell)}(x_{1}, \ldots, x_{q}) \approx v^{(\ell)}(x_{1}, \ldots, x_{q}), \quad  \ell \in \underline{k},
\]
where
\[
u^{(\ell)}(x_{1}, \ldots, x_{q})=u^{(\ell)} = \sum_{i=1}^{m} u_i^{(\ell)}(x_{1}, \ldots, x_{q}), \quad
v^{(\ell)}(x_{1}, \ldots, x_{q}) =v^{(\ell)}= \sum_{i=1}^{n} v_i^{(\ell)}(x_{1}, \ldots, x_{q}).
\]

We define \( \widehat{\mathbf{V}} \subseteq \mathbf{Sr}(2) \) to be the subvariety determined by the identities
\begin{equation}\label{22}
u^{(\ell)}(2x_{1}, \ldots, 2x_{q}) \approx v^{(\ell)}(2x_{1}, \ldots, 2x_{q}), \quad  \ell \in \underline{k}.
\end{equation}

\begin{Lemma}\label{3.1}
Let \( \mathbf{V} \in {\cal L}(\mathbf{aiSr}(2)) \). Then
\[
\widehat{\mathbf{V}} = \mathbf{V} \vee \mathbf{HSP}(Z_2, W_2, Z_7, Z_8).
\]
\end{Lemma}

\begin{Proof}
Clearly, \( \mathbf{V} \subseteq \widehat{\mathbf{V}} \), since \( \mathbf{V} \) satisfies the identities in (\ref{22}). Moreover, since \( Z_2, W_2, Z_7, Z_8 \in \widehat{\mathbf{V}} \), we have
\[
\mathbf{V} \vee \mathbf{HSP}(Z_2, W_2, Z_7, Z_8) \subseteq \widehat{\mathbf{V}}.
\]

To prove the reverse inclusion, we show that any identity \( u \approx v \) satisfied by \( \mathbf{V} \vee \mathbf{HSP}(Z_2, W_2, Z_7, Z_8) \) can be derived from the identities holding in \( \mathbf{Sr}(2) \) and those in (\ref{22}).

Let \( u = u_1 + \cdots + u_m \), \( v = v_1 + \cdots + v_n \), with \( u_i, v_j \in X^+ \). By Lemma \ref{1.1}, we only need to consider two cases:

\paragraph{Case 1.} \( m, n \geq 2 \).  
The identity \( u \approx v \) is derivable from a collection \( \Sigma = \{u^{(\ell)} \approx v^{(\ell)} \mid \ell \in \underline{k} \} \) and the defining identities of \( \mathbf{aiSr}(2) \). By \cite[Exercise II.14.11]{bur}, there exist terms \( u=t^{(1)}, t^{(2)}, \ldots, t^{(\ell)} = v \in P_f(X^+) \) such that for each \( i = 1, \ldots, \ell - 1 \), there exist \( p^{(i)}, q^{(i)}, r^{(i)} \in P_f(X^+) \), a semiring substitution \( \varphi_i \), and an identity \( {u^{(i)}}^{'} \approx {v^{(i)}}^{'} \in \Sigma \), such that
\[
t^{(i)} = p^{(i)} \varphi_i(w^{(i)}) q^{(i)} + r^{(i)}, \quad t^{(i+1)} = p^{(i)} \varphi_i(s^{(i)}) q^{(i)} + r^{(i)},
\]
where \( (w^{(i)}, s^{(i)}) = ({u^{(i)}}^{'}, {v^{(i)}}^{'}) \) or \(( {v^{(i)}}^{'}, {u^{(i)}}^{'}) \).

Define
\[
\Sigma' = \{ u(2x_1,\ldots,2x_q) \approx v(2x_1,\ldots,2x_q) \mid u(x_1,\ldots,x_q) \approx v(x_1,\ldots,x_q) \in \Sigma \}.
\]

Then, for each \( i \), we have:
\begin{align*}
t^{(i)}(2x_1,\ldots,2x_q) &= 2t^{(i)}((x_1,\ldots,x_q)) = 2(p^{(i)}\varphi_i(w^{(i)})q^{(i)} + r^{(i)}) \\
&\approx 2(p^{(i)}\varphi_i(w^{(i)})q^{(i)}) + 2r^{(i)} \\
&\approx p^{(i)}\varphi_i(w^{(i)}(2x_1,\ldots,2x_q))q^{(i)} + 2r^{(i)} \\
&\approx p^{(i)}\varphi_i(s^{(i)}(2x_1,\ldots,2x_q))q^{(i)} + 2r^{(i)} \\
&\quad (\text{since } w^{(i)}(2x_1,\ldots,2x_q) \approx s^{(i)}(2x_1,\ldots,2x_q) \in \Sigma') \\
&\approx 2(p^{(i)}\varphi_i(s^{(i)})q^{(i)} + r^{(i)}) = t^{(i+1)}(2x_1,\ldots,2x_q).
\end{align*}
By transitivity:
\[
u(2x_1,\ldots,2x_q) = t^{(1)}(2x_1,\ldots,2x_q) \approx \cdots \approx t^{(\ell)}(2x_1,\ldots,2x_q) = v(2x_1,\ldots,2x_q).
\]
Hence,
\begin{align*}
u(x_1,\ldots,x_q) &\approx 2u(x_1,\ldots,x_q)\\
&\approx u(2x_1,\ldots,2x_q)\quad (\text{by } (\ref{19}),(\ref{19'}))\\
&\approx v(2x_1,\ldots,2x_q)\\
&\approx 2v(x_1,\ldots,x_q) \quad (\text{by } (\ref{19}),(\ref{19'}))\\
&\approx v(x_1,\ldots,x_q).
\end{align*}

\paragraph{Case 2.} \( m = n = 1 \), and \( C(u) = C(v) \).  
Since \( Z_2 \models u_1 \approx v_1 \), we must have \( u_1 \neq x, v_1 \neq x \) for all \( x \in X \). Again, by \cite[Exercise II.14.11]{bur}, a derivation chain exists using identities from \( \Sigma \) and the identities determinging   \( \mathbf{Sr}(2) \). Since by Lemma \ref{2.1}, these identities are derivable from (\ref{2}), (\ref{3}), and (\ref{1}), they are also valid in \( \mathbf{Sr}(2) \), so \( u_1 \approx v_1 \) follows.

This completes the proof.
\end{Proof}

\begin{Lemma} \label{3.2}
The following equality holds:
\begin{equation} \label{24}
{\cal L}(\mathbf{Sr}(2)) = \bigcup_{\mathbf{V} \in {\cal L}(\mathbf{aiSr}(2))} [\mathbf{V}, \widehat{\mathbf{V}}].
\end{equation}
There are 64 intervals in \emph{(\ref{24})}, and each interval is a congruence class of the kernel of the complete epimorphism \( \varphi \) in \emph{(\ref{20})}.
\end{Lemma}

\begin{Proof}
First, we show that equality (\ref{24}) holds. It is straightforward to see that
\[
{\cal L}(\mathbf{Sr}(2)) = \bigcup_{\mathbf{V} \in {\cal L}(\mathbf{aiSr}(2))} \varphi^{-1}(\mathbf{V}).
\]
Thus, it suffices to show that
\begin{equation} \label{25}
\varphi^{-1}(\mathbf{V}) = [\mathbf{V}, \widehat{\mathbf{V}}],
\end{equation}
for each \( \mathbf{V} \in {\cal L}(\mathbf{aiSr}(2)) \).

Let \( \mathbf{V}_1 \) be a member of \( [\mathbf{V}, \widehat{\mathbf{V}}] \). By definition, we have
\[
\mathbf{V} \subseteq \{ E^+(S) \mid S \in \mathbf{V}_1 \} \subseteq \mathbf{V}.
\]
This implies that \( \{ E^+(S) \mid S \in \mathbf{V}_1 \} = \mathbf{V} \), so \( \varphi(\mathbf{V}_1) = \mathbf{V} \). Hence, \( \mathbf{V}_1 \in \varphi^{-1}(\mathbf{V}) \), which shows that
\[
[\mathbf{V}, \widehat{\mathbf{V}}] \subseteq \varphi^{-1}(\mathbf{V}).
\]

Conversely, if \( \mathbf{V}_1 \in \varphi^{-1}(\mathbf{V}) \), then \( \mathbf{V} = \varphi(\mathbf{V}_1) = \{ E^+(S) \mid S \in \mathbf{V}_1 \} \). Thus, \( \mathbf{V}_1 \) must belong to \( [\mathbf{V}, \widehat{\mathbf{V}}] \), showing that
\[
\varphi^{-1}(\mathbf{V}) \subseteq [\mathbf{V}, \widehat{\mathbf{V}}].
\]

Therefore, equality (\ref{25}) holds, and consequently, equality (\ref{24}) is verified.

From Lemma \ref{2.1}, we know that \( {\cal L}(\mathbf{aiSr}(2)) \) is a lattice of order 64. Hence, there are exactly 64 intervals in the expression (\ref{24}).

Next, we show that \( \varphi \) is a complete surjective morphism. First, it is easy to verify that \( \varphi \) is a complete \( \wedge \)-surjective morphism. To show that \( \varphi \) is also a complete \( \vee \)-homomorphism, consider a family of subvarieties \( (\mathbf{V}_i)_{i \in I} \) in \( {\cal L}(\mathbf{Sr}(2)) \). By the mapping (\ref{20}), for each \( i \in I \), we have
\[
\varphi(\mathbf{V}_i) \subseteq \mathbf{V}_i \subseteq \widehat{\varphi(\mathbf{V}_i)}.
\]
Furthermore, we obtain the following chain of inclusions:
\[
\bigvee_{i \in I} \varphi(\mathbf{V}_i) \subseteq \bigvee_{i \in I} \mathbf{V}_i \subseteq \bigvee_{i \in I} \widehat{\varphi(\mathbf{V}_i)} \subseteq \widehat{\bigvee_{i \in I} \varphi(\mathbf{V}_i)}.
\]
This implies that
\[
\varphi\left( \bigvee_{i \in I} \mathbf{V}_i \right) = \bigvee_{i \in I} \varphi(\mathbf{V}_i),
\]
proving that \( \varphi \) is a complete \( \vee \)-homomorphism. Therefore, \( \varphi \) is a complete surjective morphism.

Finally, from (\ref{25}), we deduce that each interval in (\ref{24}) is a congruence class of the kernel of the complete surjective morphism \( \varphi \).
\end{Proof}

To characterize the lattice of subvarieties of ${\cal L}(\mathbf{Sr}(2))$, it suffices, by Lemma \ref{3.2}, to describe the interval $[\mathbf{V}, \widehat{\mathbf{V}}]$ for each member $\mathbf{V}$ of ${\cal L}(\mathbf{Sr}(2))$. The following lemmas will aid in this task.

\begin{Lemma} \label{3.55555}  
Let $\mathbf{V}$ be a member of ${\cal L}(\mathbf{aiSr}(2))$. Then  
\[
\mathbf{V} \vee \mathbf{HSP}(Z_2, W_2, Z_7)
\]  
is the subvariety of $\widehat{\mathbf{V}}$ determined by the identity  
\begin{equation}  
2x^2 \approx 3x^2. \label{2921111}  
\end{equation}  
\end{Lemma}  

\begin{Proof}  
It is clear that both $\mathbf{V}$ and $\mathbf{HSP}(Z_2, W_2, Z_7)$ satisfy identity (\ref{2921111}). Therefore, their join $\mathbf{V} \vee \mathbf{HSP}(Z_2, W_2, Z_7)$ also satisfies (\ref{2921111}).
To complete the proof, we show that every nontrivial identity satisfied by $\mathbf{V} \vee \mathbf{HSP}(Z_2, W_2, Z_7)$ can be derived from (\ref{2921111}) together with the identities holding in $\widehat{\mathbf{V}}$.

Let $u \approx v$ be such a nontrivial identity, where  
\[
u = u_1 + u_2 + \cdots + u_m,\quad v = v_1 + v_2 + \cdots + v_n,
\]  
with $u_i, v_j \in X^+$ for $1 \leq i \leq m$, $1 \leq j \leq n$.  
Since $W_2 \models u \approx v$ but $Z_8 \not\models u \approx v$, it follows that $m, n \geq 2$.

We consider two cases based on whether $T_2$ belongs to $\mathbf{V}$.

\medskip
\noindent
\textbf{Case 1.} $T_2 \in \mathbf{V}$.  
\begin{itemize}
  \item If $|u_i| = 1$ for all $i$, then since $T_2, Z_7 \models u \approx v$, it must be that $|v_j| = 1$ for all $j$, and  
  \[
  \{ u_i \mid |u_i| = 1,\, \operatorname{occ}(u_i, u) \text{ is odd} \} = \{ v_j \mid |v_j| = 1,\, \operatorname{occ}(v_j, v) \text{ is odd} \}.
  \]  
  Therefore, $u \approx v$ is derivable from identity (\ref{162}), which holds in $\widehat{\mathbf{V}}$.

  \item If there exists some $u_i$ with $|u_i| > 1$, then by $T_2 \models u \approx v$, there exists $v_j$ such that $|v_j| > 1$.  
  By Lemma~\ref{3.1}, the identities
  \[
  u \approx u + 2u^2,\quad v \approx v + 2v^2,\quad u + 3u^2 \approx v + 3v^2
  \]
  hold in $\widehat{\mathbf{V}}$. Then, using (\ref{2921111}), we derive:
  \[
  u \approx u + 2u^2 \overset{(\ref{2921111})}{\approx} u + 3u^2 \approx v + 3v^2 \overset{(\ref{2921111})}{\approx} v + 2v^2 \approx v.
  \]
\end{itemize}

\medskip
\noindent
\textbf{Case 2.} $T_2 \not\in \mathbf{V}$.  
Then by Lemma~\ref{3.1}, the same identities as above hold in $\widehat{\mathbf{V}}$.  
Hence, the derivation proceeds identically:
\[
u \approx u + 2u^2 \overset{(\ref{2921111})}{\approx} u + 3u^2 \approx v + 3v^2 \overset{(\ref{2921111})}{\approx} v + 2v^2 \approx v.
\]

In either case, the identity $u \approx v$ is derivable from (\ref{2921111}) and the identities satisfied by $\widehat{\mathbf{V}}$. Thus, the lemma is proved.
\end{Proof}

\begin{Lemma} \label{cccc}  
Let $\mathbf{V}$ be a member of ${\cal L}(\mathbf{aiSr}(2))$. Then  
\[
\mathbf{V} \vee \mathbf{HSP}(Z_2, Z_7, Z_8)
\]  
is the subvariety of $\widehat{\mathbf{V}}$ determined by the identity  
\begin{equation}  
xy \approx 3xy. \label{c21}  
\end{equation}  
\end{Lemma}

\begin{Proof}  
Clearly, both $\mathbf{V}$ and $\mathbf{HSP}(Z_2, Z_7, Z_8)$ satisfy the identity (\ref{c21}), and thus their join $\mathbf{V} \vee \mathbf{HSP}(Z_2, Z_7, Z_8)$ also satisfies it.

To prove the lemma, it suffices to show that any nontrivial identity satisfied by $\mathbf{V} \vee \mathbf{HSP}(Z_2, Z_7, Z_8)$ can be derived from (\ref{c21}) and the identities valid in $\widehat{\mathbf{V}}$.

Let $u \approx v$ be such an identity, where  
\[
u = u_1 + u_2 + \cdots + u_m, \quad v = v_1 + v_2 + \cdots + v_n,
\]  
with $u_i, v_j \in X^+$ for $1 \leq i \leq m$, $1 \leq j \leq n$.

We consider the following cases:

\medskip
\noindent
\textbf{Case 1.} $m = n = 1$.  
Since $Z_2 \models u_1 \approx v_1$, it follows that $|u_1| \neq 1$ and $|v_1| \neq 1$.  
By Lemma~\ref{3.1}, the identity $3u_1 \approx 3v_1$ holds in $\widehat{\mathbf{V}}$.  
Hence, using (\ref{c21}), we derive:
\[
u_1 \overset{(\ref{c21})}{\approx} 3u_1 \approx 3v_1 \overset{(\ref{c21})}{\approx} v_1.
\]

\medskip
\noindent
\textbf{Case 2.} $m = 1$, $n \geq 2$.  
Since $Z_2, Z_7, Z_8 \models u_1 \approx v$, we know that:
\begin{itemize}
  \item $|u_1| \neq 1$,
  \item the number of $v_j \in v$ such that $c(v_j) \subseteq c(u_1)$ is odd,
  \item and if $|v_j| = 1$, then $\operatorname{occ}(v_j, v)$ is even.
\end{itemize}  
By Lemma~\ref{3.1}, $\widehat{\mathbf{V}}$ satisfies the identities $3u_1 \approx 3v$ and $3v \approx v$.  
Thus:
\[
u_1 \overset{(\ref{c21})}{\approx} 3u_1 \approx 3v \approx v.
\]

\medskip
\noindent
\textbf{Case 3.} $m \geq 2$, $n = 1$.  
This case is symmetric to Case 2 and follows similarly.

\medskip
\noindent
\textbf{Case 4.} $m, n \geq 2$.  
In this case, the identity $u \approx v$ is already satisfied by $\widehat{\mathbf{V}}$.

\medskip
In all cases, we conclude that the identity $u \approx v$ is derivable from (\ref{c21}) and the identities of $\widehat{\mathbf{V}}$, as required.
\end{Proof}

\begin{Lemma} \label{3.64321}  
Let $\mathbf{V}$ be a member of ${\cal L}(\mathbf{aiSr}(2))$. Then  
\[
\mathbf{V} \vee \mathbf{HSP}(Z_2, W_2)
\]  
is the subvariety of $\widehat{\mathbf{V}}$ determined by the identity  
\begin{equation}
x + u' \approx 2x + v', \quad \text{where } \mathbf{V} \models u' \approx v',\, x \notin u',\, x \notin v'. \label{292121}
\end{equation}  
\end{Lemma}

\begin{Proof}  
It is straightforward to verify that both $\mathbf{V}$ and $\mathbf{HSP}(Z_2, W_2)$ satisfy the identity (\ref{292121}), and hence so does their join $\mathbf{V} \vee \mathbf{HSP}(Z_2, W_2)$.

To complete the proof, we show that every nontrivial identity satisfied by $\mathbf{V} \vee \mathbf{HSP}(Z_2, W_2)$ is derivable from (\ref{292121}) and the identities satisfied by $\widehat{\mathbf{V}}$.

Let $u \approx v$ be such an identity, where  
\[
u = u_1 + u_2 + \cdots + u_m, \quad v = v_1 + v_2 + \cdots + v_n, \quad u_i, v_j \in X^+, \ 1 \leq i \leq m, \ 1 \leq j \leq n.
\]  
Since $W_2 \models u \approx v$, we consider the following two main cases:

\medskip
\noindent
\textbf{Case 1.} $m = n = 1$ and $c(u_1) = c(v_1)$.  
Since $Z_2 \models u_1 \approx v_1$, $|u_1|, |v_1| \geq 2$. It follows that $\widehat{\mathbf{V}} \models u_1 \approx v_1$.

\medskip
\noindent
\textbf{Case 2.} $m, n \geq 2$. We divide this into two subcases:

\begin{itemize}
  \item[\textbf{(a)}] $T_2 \in \mathbf{V}$ and there exists some $u_i$ with $|u_i| \geq 2$, or $T_2 \notin \mathbf{V}$.  
  By Lemma~\ref{3.1}, $\widehat{\mathbf{V}}$ satisfies:
  \[
  u \approx u + u^2, \quad v \approx v + v^2, \quad 2u + u^2 \approx 2v + v^2.
  \]  
  Applying (\ref{292121}), we obtain:
  \[
  u \approx u + u^2 \overset{(\ref{292121})}{\approx} 2u + u^2 \approx 2v + v^2 \overset{(\ref{292121})}{\approx} v + v^2 \approx v.
  \]

  \item[\textbf{(b)}] $T_2 \in \mathbf{V}$ and $|u_i| = 1$ for all $i$.  
  Since either $Z_7 \not\models u \approx v$ or $Z_8 \not\models u \approx v$, we may, without loss of generality, assume that  
  \[
  u = x + u'', \quad v = 2x + v'', \quad \text{with } x \notin u'',\, x \notin v''.
  \]  
  Then the identity $u \approx v$ is an instance of (\ref{292121}), and thus holds.
\end{itemize}

In all cases, the identity $u \approx v$ is derivable from (\ref{292121}) and identities satisfied by $\widehat{\mathbf{V}}$. This completes the proof.
\end{Proof}

\begin{Lemma} \label{ccccc}  
Let $\mathbf{V}$ be a member of ${\cal L}(\mathbf{aiSr}(2))$. Then  
\[
\mathbf{V} \vee \mathbf{HSP}(Z_2, Z_7)
\]  
is the subvariety of $\widehat{\mathbf{V}}$ determined by the identity  
\begin{equation}
xy \approx 2xy. \label{c21'}
\end{equation}  
\end{Lemma}

\begin{Proof}  
It is clear that both $\mathbf{V}$ and $\mathbf{HSP}(Z_2, Z_7)$ satisfy the identity (\ref{c21'}), and thus their join $\mathbf{V} \vee \mathbf{HSP}(Z_2, Z_7)$ does as well.

To complete the proof, it suffices to show that every nontrivial identity satisfied by $\mathbf{V} \vee \mathbf{HSP}(Z_2, Z_7)$ can be derived from (\ref{c21'}) and the identities valid in $\widehat{\mathbf{V}}$.

Let $u \approx v$ be such an identity, where  
\[
u = u_1 + u_2 + \cdots + u_m,\quad v = v_1 + v_2 + \cdots + v_n,\quad u_i, v_j \in X^+,\ 1 \leq i \leq m,\ 1 \leq j \leq n.
\]

We consider the following cases:

\medskip
\noindent
\textbf{Case 1.} $m = n = 1$.  
Since $Z_2 \models u_1 \approx v_1$, we have $|u_1|, |v_1| \geq 2$. By Lemma~\ref{3.1}, $\widehat{\mathbf{V}} \models 2u_1 \approx 2v_1$.  
Using identity (\ref{c21'}), we obtain:
\[
u_1 \overset{(\ref{c21'})}{\approx} 2u_1 \approx 2v_1 \overset{(\ref{c21'})}{\approx} v_1.
\]

\medskip
\noindent
\textbf{Case 2.} $m = 1$, $n \geq 2$.  
Since $Z_2, Z_7 \models u_1 \approx v$, we have:
\begin{itemize}
  \item $|u_1| \geq 2$,
  \item for any $v_j \in v$ with $|v_j| = 1$, the number of occurrences $\operatorname{occ}(v_j, v)$ is even.
\end{itemize}
By Lemma~\ref{3.1}, $\widehat{\mathbf{V}} \models 2u_1 \approx 2v$.  
Then:
\[
u_1 \overset{(\ref{c21'})}{\approx} 2u_1 \approx 2v \overset{(\ref{162}),(\ref{c21'})}{\approx} v.
\]

\medskip
\noindent
\textbf{Case 3.} $m \geq 2$, $n = 1$.  
This case is symmetric to Case 2 and follows similarly.

\medskip
\noindent
\textbf{Case 4.} $m, n \geq 2$.  
By Lemma~\ref{3.1}, $\widehat{\mathbf{V}}$ satisfies $2u \approx 2v$.  
Applying identity (\ref{c21'}), we obtain:
\[
u \overset{(\ref{c21'})}{\approx} 2u \approx 2v \overset{(\ref{c21'})}{\approx} v.
\]

\medskip
In all cases, the identity $u \approx v$ is derivable from (\ref{c21'}) and the identities holding in $\widehat{\mathbf{V}}$, completing the proof.
\end{Proof}

\begin{Lemma} \label{3.35672}  
Let $\mathbf{V}$ be a member of ${\cal L}(\mathbf{aiSr}(2))$ such that $N_2 \notin \mathbf{V}$. Then  
\[
\mathbf{V} \vee \mathbf{HSP}(Z_2, Z_8)
\]  
is the subvariety of $\widehat{\mathbf{V}}$ determined by the identity  
\begin{equation}
x^2 \approx x + 2x^2. \label{27,,}
\end{equation}  
\end{Lemma}

\begin{Proof}  
Clearly, both $\mathbf{V}$ and $\mathbf{HSP}(Z_2, Z_8)$ satisfy the identity (\ref{27,,}), and hence so does their join $\mathbf{V} \vee \mathbf{HSP}(Z_2, Z_8)$.

To prove the lemma, it suffices to show that every nontrivial identity satisfied by $\mathbf{V} \vee \mathbf{HSP}(Z_2, Z_8)$ is derivable from (\ref{27,,}) together with the identities holding in $\widehat{\mathbf{V}}$.

Let $u \approx v$ be such an identity, where  
\[
u = u_1 + u_2 + \cdots + u_m, \quad v = v_1 + v_2 + \cdots + v_n, \quad u_i, v_j \in X^+, \ 1 \leq i \leq m, \ 1 \leq j \leq n.
\]  
Since $Z_2 \models u \approx v$, we must have $u \neq x$ and $v \neq x$ for any $x \in X$. We consider the following cases:

\medskip
\noindent
\textbf{Case 1.} $m = n = 1$, with $|u_1|, |v_1| \geq 2$.  
By Lemma~\ref{3.1}, $\widehat{\mathbf{V}}$ satisfies:
\[
u_1 \approx u_1^2,\quad v_1 \approx v_1^2,\quad u_1 + 2u_1^2 \approx v_1 + 2v_1^2.
\]
Applying (\ref{27,,}), we derive:
\[
u_1 \approx u_1^2 \overset{(\ref{27,,})}{\approx} u_1 + 2u_1^2 \approx v_1 + 2v_1^2 \overset{(\ref{27,,})}{\approx} v_1^2 \approx v_1.
\]

\medskip
\noindent
\textbf{Case 2.} $m = 1$, $n \geq 2$, with $|u_1| \geq 2$.  
By Lemma~\ref{3.1}, $\widehat{\mathbf{V}}$ satisfies:
\[
u_1 \approx u_1^2,\quad v \approx v + 2v^2,\quad u_1 + 2u_1^2 \approx v^2.
\]
Then:
\[
u_1 \approx u_1^2 \overset{(\ref{27,,})}{\approx} u_1 + 2u_1^2 \approx v^2 \overset{(\ref{27,,})}{\approx} v + 2v^2 \approx v.
\]

\medskip
\noindent
\textbf{Case 3.} $n = 1$, $m \geq 2$, and $|v_1| \geq 2$.  
This case is symmetric to Case 2 and follows the same reasoning.

\medskip
\noindent
\textbf{Case 4.} $m, n \geq 2$. We consider two subcases:

\smallskip
\noindent
\textbf{Subcase 4.1.} There exists $u_i$ with $|u_i| \geq 2$.  
Then by Lemma~\ref{3.1}, $\widehat{\mathbf{V}}$ satisfies:
\[
u \approx u + 2u^2,\quad v \approx v + 2v^2,\quad u^2 \approx v^2.
\]
Thus:
\[
u \approx u + 2u^2 \overset{(\ref{27,,})}{\approx} u^2 \approx v^2 \overset{(\ref{27,,})}{\approx} v + 2v^2 \approx v.
\]

\smallskip
\noindent
\textbf{Subcase 4.2.} All $u_i$ satisfy $|u_i| = 1$.  
Note that $Z_8 \models u \approx v$. We distinguish two possibilities:
\begin{itemize}
  \item If $\mathbf{V} \cap \{L_2, R_2, M_2, D_2\} \neq \varnothing$, then $u \approx v$ follows by identity (\ref{162}).
  \item If $\mathbf{V} \cap \{L_2, R_2, M_2, D_2\} = \varnothing$, then $u \approx v$ is already satisfied in $\widehat{\mathbf{V}}$.
\end{itemize}

\medskip
In all cases, the identity $u \approx v$ is derivable from (\ref{27,,}) and the identities of $\widehat{\mathbf{V}}$, completing the proof.
\end{Proof}

\begin{Lemma} \label{3.4}  
Let $\mathbf{V}$ be a member of ${\cal L}(\mathbf{HSP}(L_2,R_2,M_2,D_2))$. Then  
\[
\mathbf{V} \vee \mathbf{HSP}(W_2, Z_8)
\]  
is the subvariety of $\widehat{\mathbf{V}}$ determined by the identity  
\begin{equation}
x^2 \approx x. \label{29}
\end{equation}  
\end{Lemma}

\begin{Proof}  
It is easy to verify that both $\mathbf{V}$ and $\mathbf{HSP}(W_2, Z_8)$ satisfy the identity (\ref{29}), and hence their join $\mathbf{V} \vee \mathbf{HSP}(W_2, Z_8)$ does as well.

To complete the proof, it suffices to show that every nontrivial identity satisfied by $\mathbf{V} \vee \mathbf{HSP}(W_2, Z_8)$ is derivable from (\ref{29}) and the identities valid in $\widehat{\mathbf{V}}$.

Let $u \approx v$ be such an identity, where  
\[
u = u_1 + u_2 + \cdots + u_m, \quad v = v_1 + v_2 + \cdots + v_n, \quad u_i, v_j \in X^+, \ 1 \leq i \leq m,\ 1 \leq j \leq n.
\]

We analyze the following cases:

\medskip
\noindent
\textbf{Case 1.} $m = n = 1$, and $c(u_1) = c(v_1)$.  
By Lemma~\ref{3.1}, we know that $\widehat{\mathbf{V}} \models u_1^2 \approx v_1^2$.  
Then, by identity (\ref{29}), we have:  
\[
u_1 \overset{(\ref{29})}{\approx} u_1^2 \approx v_1^2 \overset{(\ref{29})}{\approx} v_1.
\]

\medskip
\noindent
\textbf{Case 2.} $m, n \geq 2$.  
By Lemma~\ref{3.1}, $\widehat{\mathbf{V}}$ satisfies the identity $u^2 \approx v^2$.  
Then:
\[
u \overset{(\ref{29})}{\approx} u^2 \approx v^2 \overset{(\ref{29})}{\approx} v.
\]

\medskip
In both cases, the identity $u \approx v$ is derivable from (\ref{29}) and the identities valid in $\widehat{\mathbf{V}}$, completing the proof.
\end{Proof}

\begin{Lemma} \label{3.555}  
Let $\mathbf{V}$ be a member of ${\cal L}(\mathbf{aiSr}(2))$ such that $T_2 \notin \mathbf{V}$. Then  
\[
\mathbf{V} \vee \mathbf{HSP}(Z_7, Z_8)
\]  
is the subvariety of $\widehat{\mathbf{V}}$ determined by the identity  
\begin{equation}
x \approx x + 2x^2. \label{29211}
\end{equation}  
\end{Lemma}

\begin{Proof}  
It is straightforward to verify that both $\mathbf{V}$ and $\mathbf{HSP}(Z_7, Z_8)$ satisfy the identity (\ref{29211}), and hence their join $\mathbf{V} \vee \mathbf{HSP}(Z_7, Z_8)$ also satisfies it.

We aim to show that every nontrivial identity satisfied by $\mathbf{V} \vee \mathbf{HSP}(Z_7, Z_8)$ can be derived from (\ref{29211}) and the identities valid in $\widehat{\mathbf{V}}$.

Let $u \approx v$ be such an identity, where  
\[
u = u_1 + u_2 + \cdots + u_m, \quad v = v_1 + v_2 + \cdots + v_n, \quad u_i, v_j \in X^+, \ 1 \leq i \leq m,\ 1 \leq j \leq n.
\]  
By Lemma~\ref{3.1},  $\widehat{\mathbf{V}}$ satisfies the identity  
\[
u + 2u^2 \approx v + 2v^2.
\]  
Then, using (\ref{29211}), we obtain:
\[
u \overset{(\ref{29211})}{\approx} u + 2u^2 \approx v + 2v^2 \overset{(\ref{29211})}{\approx} v.
\]

This completes the proof.
\end{Proof}

\begin{Lemma} \label{3.3567}
Let $\mathbf{V}$ be a member of ${\cal L}(\mathbf{aiSr}(2))$ such that $N_2 \notin \mathbf{V}$. Then
\[
\mathbf{V} \vee \mathbf{HSP}(Z_2)
\]
is the subvariety of $\widehat{\mathbf{V}}$ determined by the identity
\begin{equation}
x^2 \approx x + x^2. \label{27}
\end{equation}
\end{Lemma}

\begin{Proof}
It is clear that both $\mathbf{V}$ and $\mathbf{HSP}(Z_2)$ satisfy the identity (\ref{27}), and hence so does $\mathbf{V} \vee \mathbf{HSP}(Z_2)$.

We now show that every nontrivial identity satisfied by $\mathbf{V} \vee \mathbf{HSP}(Z_2)$ can be derived from (\ref{27}) together with the identities valid in $\widehat{\mathbf{V}}$.

Let $u \approx v$ be such an identity, where
\[
u = u_1 + u_2 + \cdots + u_m,  v = v_1 + v_2 + \cdots + v_n, u_i, v_j \in X^+, \ 1 \leq i \leq m,\ 1 \leq j \leq n.
\]
Since $Z_2 \models u \approx v$, we have $u \neq x$ and $v \neq x$ for any $x \in X$.
We consider the following cases:

\textbf{Case 1.} $\forall u_i \in u$, $|u_i| = 1$. We consider the following subcases:

\begin{itemize}
  \item[$\bullet$] If \(T_2\notin \mathbf{V} \), then by Lemma~\ref{3.1}, \( \widehat{\mathbf{V}} \) satisfies the identities
  \[
  u \approx u + 2u^2, v \approx v + 2v^2,  2u^2 \approx 2v^2.
  \]
  Then,
  \[
  u \approx u + 2u^2 \overset{(\ref{27})}{\approx} 2u^2 \approx 2v^2 \overset{(\ref{27})}{\approx} v + 2v^2 \approx v.
  \] 

  \item[$\bullet$]  If \( T_2\in \mathbf{V} \), then by Lemma~\ref{3.1}, \( \widehat{\mathbf{V}} \) satisfies the identities
  \[
  u \approx u + 2v, v \approx v + 2u.
  \]
  Then,
  \[
  u \approx u + 2v\approx  v + 2u \approx v.
  \] 
\end{itemize}

\textbf{Case 2.} $\exists u_i \in u$ such that $|u_i| \geq 2$. We discuss the following subcases:

\begin{itemize}
  \item[$\bullet$] $m, n \geq 2$.  
  By Lemma~\ref{3.1}, $\widehat{\mathbf{V}}$ satisfies the identities $u \approx u + 2u^2$, $v \approx v + 2v^2$, and $2u^2 \approx 2v^2$.  
  Then,
  \[
  u \approx u + 2u^2 \overset{(\ref{27})}{\approx} 2u^2 \approx 2v^2 \overset{(\ref{27})}{\approx} v + 2v^2 \approx v.
  \]

  \item[$\bullet$] $m = n = 1$.  
  By Lemma~\ref{3.1}, $\widehat{\mathbf{V}}$ satisfies $u_1 \approx u_1^2$, $v_1 \approx v_1^2$, and $u_1 + u_1^2 \approx v_1 + v_1^2$.  
  Thus,
  \[
  u_1 \approx u_1^2 \overset{(\ref{27})}{\approx} u_1 + u_1^2 \approx v_1 + v_1^2 \overset{(\ref{27})}{\approx} v_1^2 \approx v_1.
  \]

  \item[$\bullet$] $m = 1$, $n > 1$.  
  Then $\widehat{\mathbf{V}}$ satisfies $u_1 \approx u_1^2$, $v \approx v + 2v^2$, and $2u_1^2 \approx 2v^2$. Hence,
  \[
  u_1 \approx u_1^2 \overset{(\ref{27})}{\approx} u_1 + u_1^2 \approx 2u_1^2 \approx 2v^2 \overset{(\ref{27})}{\approx} v + 2v^2 \approx v.
  \]

  \item[$\bullet$] $m > 1$, $n = 1$.  
  This case is symmetric to the previous one and follows analogously.
\end{itemize}
\end{Proof}

\begin{Lemma} \label{3.55}
Let $\mathbf{V}$ be a member of ${\cal L}(\mathbf{aiSr}(2))$ such that $T_2 \in \mathbf{V}$. Then
\[
\mathbf{V} \vee \mathbf{HSP}(Z_7)
\]
is the subvariety of $\widehat{\mathbf{V}}$ determined by the identity
\begin{equation}
x \approx x + x^2. \label{2921}
\end{equation}
\end{Lemma}

\begin{Proof}
It is easy to verify that both $\mathbf{V}$ and $\mathbf{HSP}(Z_7)$ satisfy the identity (\ref{2921}). Consequently, their join $\mathbf{V} \vee \mathbf{HSP}(Z_7)$ also satisfies (\ref{2921}).

We now prove that every nontrivial identity satisfied by $\mathbf{V} \vee \mathbf{HSP}(Z_7)$ is derivable from (\ref{2921}) and the identities valid in $\widehat{\mathbf{V}}$.

Let $u \approx v$ be such an identity, where
\[
u = u_1 + u_2 + \cdots + u_m, \quad v = v_1 + v_2 + \cdots + v_n, \quad u_i, v_j \in X^+, \ 1 \leq i \leq m,\ 1 \leq j \leq n.
\]
By Lemma~\ref{3.1}, the identity
\[
u + u^2 \approx v + v^2
\]
holds in $\widehat{\mathbf{V}}$.

Applying (\ref{2921}), we derive:
\[
u \overset{(\ref{2921})}{\approx} u + u^2 \approx v + v^2 \overset{(\ref{2921})}{\approx} v.
\]
This completes the proof.
\end{Proof}

\begin{Lemma} \label{3.5555}
Let $\mathbf{V}$ be a member of ${\cal L}(\mathbf{aiSr}(2))$ and $T_2 \notin \mathbf{V}$. Then
\[
\mathbf{V} \vee \mathbf{HSP}(Z_8)
\]
is the subvariety of $\widehat{\mathbf{V}}$ determined by the identity
\begin{equation}
x \approx 2x + x^2. \label{292111}
\end{equation}
\end{Lemma}

\begin{Proof}
It is easy to verify that both $\mathbf{V}$ and $\mathbf{HSP}(Z_8)$ satisfy the identity (\ref{292111}). Consequently, their join $\mathbf{V} \vee \mathbf{HSP}(Z_8)$ also satisfies (\ref{292111}).

We now prove that every nontrivial identity satisfied by $\mathbf{V} \vee \mathbf{HSP}(Z_8)$ is derivable from (\ref{292111}) and the identities valid in $\widehat{\mathbf{V}}$.

Let $u \approx v$ be such an identity, where
\[
u = u_1 + u_2 + \cdots + u_m, \quad v = v_1 + v_2 + \cdots + v_n, \quad u_i, v_j \in X^+, \ 1 \leq i \leq m,\ 1 \leq j \leq n.
\]
By Lemma~\ref{3.1}, the identity
\[
2u + u^2 \approx 2v + v^2
\]
holds in $\widehat{\mathbf{V}}$.

Applying (\ref{292111}), we derive:
\[
u \overset{(\ref{292111})}{\approx} 2u + u^2 \approx 2v + v^2 \overset{(\ref{292111})}{\approx} v.
\]
This completes the proof.
\end{Proof}

\begin{Theorem} \label{3.5} 
Let $\mathbf{V}$ be an element of ${\cal L}(\mathbf{aiSr}(2))$. The interval $[\mathbf{V}, \widehat{\mathbf{V}}]$ in ${\cal L}(\mathbf{Sr}(2))$ is depicted in Fig. 1. Specifically:
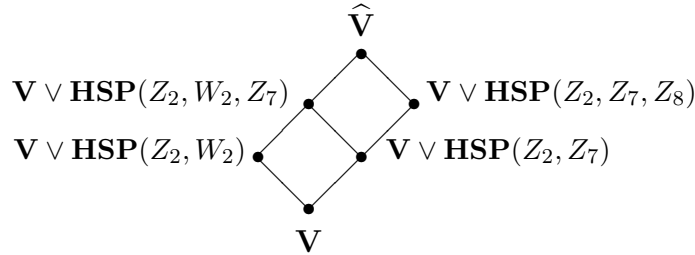
\[
\setlength{\unitlength}{1.4cm}
\begin{picture}(12,3.5)
\put(6,1.4){\line(-1,1){0.5}}
\put(5.5,0.9){\line(1,1){1}}
\put(5.5,0.9){\line(-1,1){0.5}}
\put(5,1.4){\line(1,1){1}}
\put(6.5,1.9){\line(-1,1){0.5}}
\multiput(6,2.375)(1,1){1}{\circle*{0.1}}
\multiput(5.02,1.4)(1,1){1}{\circle*{0.1}}
\multiput(6,1.4)(1,1){1}{\circle*{0.1}}
\multiput(5.5,0.9)(1,1){1}{\circle*{0.1}}
\multiput(5.5,1.9)(1,1){1}{\circle*{0.1}}
\multiput(6.5,1.9)(1,1){1}{\circle*{0.1}}
\put(5.5,0.6){\makebox(0,0){$\mathbf{V}$}}
\put(6.0,2.7){\makebox(0,0){$\mathbf{\widehat{V}}$}}
\put(3.8,1.45){\makebox(0,0){\small{$\mathbf{V} \vee \mathbf{HSP}(Z_2, W_2)$}}}
\put(7.3,1.45){\makebox(0,0){\small{$\mathbf{V} \vee \mathbf{HSP}(Z_2, Z_7)$}}}
\put(7.9,2){\makebox(0,0){\small{$\mathbf{V} \vee \mathbf{HSP}(Z_2, Z_7, Z_8)$}}}
\put(4,2){\makebox(0,0){\small{$\mathbf{V} \vee \mathbf{HSP}(Z_2, W_2, Z_7)$}}}
\put(5.5,0){\makebox(0,0){\text{Case 1:} $N_2, T_2 \in \mathbf{V}$}}
\end{picture}
\]
\[
\setlength{\unitlength}{1.4cm}
\begin{picture}(12,3.5)
\put(6,1.4){\line(-1,1){0.5}}
\put(6,1.4){\line(1,1){0.5}}
\put(6,-0.4){\line(0,1){2.8}}
\put(6,0.4){\line(-1,1){1}}
\put(6,0.4){\line(1,1){1}}
\put(5,1.4){\line(1,1){1}}
\put(7,1.4){\line(-1,1){1}}
\multiput(6,2.375)(1,1){1}{\circle*{0.1}}
\multiput(5.02,1.4)(1,1){1}{\circle*{0.1}}
\multiput(6,1.4)(1,1){1}{\circle*{0.1}}
\multiput(6.98,1.4)(1,1){1}{\circle*{0.1}}
\multiput(6.0,0.42)(1,1){1}{\circle*{0.1}}
\multiput(6.0,-0.38)(1,1){1}{\circle*{0.1}}
\multiput(5.5,1.9)(1,1){1}{\circle*{0.1}}
\multiput(6.5,1.9)(1,1){1}{\circle*{0.1}}
\put(7.1,0.3){\makebox(0,0){$\mathbf{V} \vee \mathbf{HSP}(Z_2)$}}
\put(6.0,-0.7){\makebox(0,0){$\mathbf{V}$}}
\put(6.0,2.7){\makebox(0,0){$\mathbf{\widehat{V}}$}}
\put(3.8,1.45){\makebox(0,0){\small{$\mathbf{V} \vee \mathbf{HSP}(Z_2, W_2)$}}}
\put(8.2,1.45){\makebox(0,0){\small{$\mathbf{V} \vee \mathbf{HSP}(Z_2, Z_8)$}}}
\put(7.9,2){\makebox(0,0){\small{$\mathbf{V} \vee \mathbf{HSP}(Z_2, Z_7, Z_8)$}}}
\put(4,2){\makebox(0,0){\small{$\mathbf{V} \vee \mathbf{HSP}(Z_2, W_2, Z_7)$}}}
\put(5.5,-1){\makebox(0,0){\text{Case 2:} $N_2 \notin \mathbf{V},T_2 \in \mathbf{V}$}}
\end{picture}
\]

\[
\setlength{\unitlength}{1.4cm}
\begin{picture}(12,3.5)
\put(7,0.4){\line(-2,1){2}}
\put(5.5,1.9){\line(1,-1){1.5}}
\put(6,1.4){\line(1,1){0.5}}
\put(6.5,0.9){\line(1,1){0.5}}
\put(5,1.4){\line(1,1){1}}
\put(6,2.4){\line(1,-1){1.5}}
\put(7.5,0.9){\line(-1,-1){0.5}}
\multiput(6,2.375)(1,1){1}{\circle*{0.1}}
\multiput(5.02,1.4)(1,1){1}{\circle*{0.1}}
\multiput(6,1.4)(1,1){1}{\circle*{0.1}}
\multiput(6.98,1.4)(1,1){1}{\circle*{0.1}}
\multiput(7.0,0.42)(1,1){1}{\circle*{0.1}}
\multiput(6.5,0.9)(1,1){1}{\circle*{0.1}}
\multiput(7.5,0.9)(1,1){1}{\circle*{0.1}}
\multiput(5.5,1.9)(1,1){1}{\circle*{0.1}}
\multiput(6.5,1.9)(1,1){1}{\circle*{0.1}}
\put(7,0.1){\makebox(0,0){$\mathbf{V}$}}
\put(8.5,0.9){\makebox(0,0){$\mathbf{V} \vee \mathbf{HSP}(Z_8)$}}
\put(5.5,0.7){\makebox(0,0){$\mathbf{V} \vee \mathbf{HSP}(Z_7)$}}
\put(6.0,2.7){\makebox(0,0){$\mathbf{\widehat{V}}$}}
\put(3.8,1.45){\makebox(0,0){\small{$\mathbf{V} \vee \mathbf{HSP}(Z_2, W_2)$}}}
\put(8.2,1.45){\makebox(0,0){\small{$\mathbf{V} \vee \mathbf{HSP}(Z_7, Z_8)$}}}
\put(7.9,2){\makebox(0,0){\small{$\mathbf{V} \vee \mathbf{HSP}(Z_2, Z_7, Z_8)$}}}
\put(4,2){\makebox(0,0){\small{$\mathbf{V} \vee \mathbf{HSP}(Z_2, W_2, Z_7)$}}}
\put(5.5,-0.4){\makebox(0,0){\text{Case 3:} $N_2 \in \mathbf{V},T_2 \notin \mathbf{V}$}}
\end{picture}
\]
\[
\setlength{\unitlength}{1.4cm}
\begin{picture}(12,3.5)
\put(6,1.4){\line(-1,1){0.5}}
\put(6,1.4){\line(1,1){0.5}}
\put(6,0.4){\line(0,1){2}}
\put(6,0.4){\line(-1,1){1}}
\put(6,0.4){\line(1,1){1}}
\put(5,1.4){\line(1,1){1}}
\put(7,1.4){\line(-1,1){1}}
\multiput(6,2.375)(1,1){1}{\circle*{0.1}}
\multiput(5.02,1.4)(1,1){1}{\circle*{0.1}}
\multiput(6,1.4)(1,1){1}{\circle*{0.1}}
\multiput(6.98,1.4)(1,1){1}{\circle*{0.1}}
\multiput(6.0,0.42)(1,1){1}{\circle*{0.1}}
\multiput(5.5,1.9)(1,1){1}{\circle*{0.1}}
\multiput(6.5,1.9)(1,1){1}{\circle*{0.1}}
\put(6.0,0.1){\makebox(0,0){$\mathbf{V}$}}
\put(6.0,2.7){\makebox(0,0){$\mathbf{\widehat{V}}$}}
\put(3.8,1.45){\makebox(0,0){\small{$\mathbf{V} \vee \mathbf{HSP}(Z_2, W_2)$}}}
\put(8.2,1.45){\makebox(0,0){\small{$\mathbf{V} \vee \mathbf{HSP}(W_2, Z_8)$}}}
\put(7.9,2){\makebox(0,0){\small{$\mathbf{V} \vee \mathbf{HSP}(Z_2, Z_7, Z_8)$}}}
\put(4,2){\makebox(0,0){\small{$\mathbf{V} \vee \mathbf{HSP}(Z_2, W_2, Z_7)$}}}
\put(5.5,-0.4){\makebox(0,0){\text{Case 4:} $N_2, T_2 \notin \mathbf{V}$}}
\end{picture}
\]
\vskip 30pt
\begin{center}
\textbf{Fig. 1: The interval $[\mathbf{V}, \widehat{\mathbf{V}}]$}
\end{center}

\end{Theorem}

\begin{Proof}
Suppose that $\mathbf{V}_1 \in [\mathbf{V}, \widehat{\mathbf{V}}]$ with $\mathbf{V}_1 \ne \mathbf{V}$ and $\mathbf{V}_1 \ne \widehat{\mathbf{V}}$. Then there exists a nontrivial identity $u \approx v$ such that $u \approx v$ holds in $\mathbf{V}_1$ but does not hold in $\widehat{\mathbf{V}}$. Moreover, we are given that $\mathbf{V}_1$ does not satisfy the identity $2x \approx x$. 

By Lemma~\ref{3.1}, it suffices to consider the following cases.

\textbf{Case 1.} Suppose that $Z_2, W_2, Z_7 \models u \approx v$ but $Z_8 \not\models u \approx v$. Then the identity $u \approx v$ satisfies one of the following two conditions:
\begin{itemize}
    \item $m, n \geq 2$, there exists some $u_i \in u$ such that $|u_i| = 1$, and both $\operatorname{occ}(u_i, u)$ and $\operatorname{occ}(u_i, v)$ are odd. Moreover, there exists a nonempty subset $A \subseteq X$ such that the cardinalities
    \[
        \left|\left\{ i \mid u_i \in u,\ c(u_i) \subseteq A \right\}\right| \quad \text{and} \quad \left|\left\{ j \mid v_j \in v,\ c(v_j) \subseteq A \right\}\right|
    \]
    are of different parity (i.e., one is odd, the other is even);

    \item $m, n \geq 2$, and for all $u_i \in u$, we have $|u_i| > 1$. Furthermore, there exists a nonempty subset $A \subseteq X$ such that the cardinalities
    \[
        \left|\left\{ i \mid u_i \in u,\ c(u_i) \subseteq A \right\}\right| \quad \text{and} \quad \left|\left\{ j \mid v_j \in v,\ c(v_j) \subseteq A \right\}\right|
    \]
    are again of different parity.
\end{itemize}

In either case, the identity $u \approx v$ implies the identity \(2x^2 \approx 3x^2\). By Lemma~\ref{3.55555}, it follows that \(\mathbf{V}_1\) is a subvariety of \(\mathbf{V} \vee \mathbf{HSP}(Z_2, W_2, Z_7)\).

On the other hand, since \(\mathbf{V}_1 \models 2x^2 \approx 3x^2\) but \(\mathbf{V}_1 \not\models x \approx 2x\), we conclude that \(Z_2, W_2, Z_7 \in \mathbf{V}_1\), and thus \(\mathbf{V} \vee \mathbf{HSP}(Z_2, W_2, Z_7) \subseteq \mathbf{V}_1\).

Therefore, we conclude that:
\[
\mathbf{V}_1 = \mathbf{V} \vee \mathbf{HSP}(Z_2, W_2, Z_7).
\]

\textbf{Case 2.} Suppose that $Z_2, W_2, Z_8 \models u \approx v$ but $Z_7 \not\models u \approx v$. Since $Z_2$ and $W_2$ satisfy the identity $u \approx v$ while $Z_7$ does not, the identity $u \approx v$ must fall into one of the following two cases:
\begin{itemize}
    \item $m, n \geq 2$, there exists some $u_i \in u$ with $|u_i| = 1$ such that $\operatorname{occ}(u_i, u)$ is odd and $\operatorname{occ}(u_i, v)$ is even;
    \item $m, n \geq 2$, there exists some $v_j \in v$ with $|v_j| = 1$ such that $\operatorname{occ}(v_j, v)$ is odd and $\operatorname{occ}(v_j, u)$ is even.
\end{itemize}

However, since $Z_8 \models u \approx v$, it follows that for every nonempty subset $A \subseteq X$, the number of occurrences of blocks whose content is contained in $A$ is odd in $u$ if and only if it is odd in $v$, that is,
\[
\left| \left\{ i \mid u_i \in u,\ c(u_i) \subseteq A \right\} \right| \text{ is odd} \iff \left| \left\{ j \mid v_j \in v,\ c(v_j) \subseteq A \right\} \right| \text{ is odd}.
\]
Now take \( A = c(u_i) \) or \( A = c(v_j) \) depending on which of the two cases above applies. This leads to a contradiction with the parity condition derived from the failure in \( Z_7 \).

Therefore, we conclude that:
\[
\mathbf{V} \vee \mathbf{HSP}(Z_2, W_2, Z_8) = \widehat{\mathbf{V}}.
\]

\textbf{Case 3.} Suppose that $Z_2, Z_7, Z_8 \models u \approx v$ but $W_2 \not\models u \approx v$. Then the identity $u \approx v$ must satisfy one of the following two conditions:
\begin{itemize}
    \item $m = 1$, $n \geq 2$, and $|u_1| \geq 2$. Moreover, both
    \[
        \left| \left\{ j \mid v_j \in v,\ c(v_j) \subseteq c(u_1) \right\} \right|
        \quad \text{and} \quad
        \left| \left\{ j \mid v_j \in v \right\} \right|
    \]
    are odd, and for any $v_j \in v$ with $|v_j| = 1$, we have $\operatorname{occ}(v_j, v)$ even;

    \item $n = 1$, $m \geq 2$, and $|v_1| \geq 2$. Moreover, both
    \[
        \left| \left\{ i \mid u_i \in u,\ c(u_i) \subseteq c(v_1) \right\} \right|
        \quad \text{and} \quad
        \left| \left\{ i \mid u_i \in u \right\} \right|
    \]
    are odd, and for any $u_j \in u$ with $|u_j| = 1$, we have $\operatorname{occ}(u_j, u)$ even.
\end{itemize}

In either of the above cases, the identity $u \approx v$ implies the identity \(xy \approx 3xy\). By Lemma~\ref{cccc}, it follows that \(\mathbf{V}_1\) is a subvariety of \(\mathbf{V} \vee \mathbf{HSP}(Z_2, Z_7, Z_8)\).

On the other hand, since \(\mathbf{V}_1 \models xy \approx 3xy\) and \(\mathbf{V}_1 \not\models x \approx 2x\), we conclude that \(Z_2, Z_7, Z_8 \in \mathbf{V}_1\), and hence \(\mathbf{V} \vee \mathbf{HSP}(Z_2, Z_7, Z_8) \subseteq \mathbf{V}_1\).

Therefore, we obtain:
\[
\mathbf{V}_1 = \mathbf{V} \vee \mathbf{HSP}(Z_2, Z_7, Z_8).
\]

\textbf{Case 4.} Suppose that $W_2, Z_7, Z_8 \models u \approx v$ but $Z_2 \not\models u \approx v$. Since $Z_2 \not\models u \approx v$, it must be that $u = x$ or $v = x$.

On the other hand, since $W_2, Z_7, Z_8 \models u \approx v$, we must have $u = v = x$, which contradicts the assumption that $Z_2 \not\models u \approx v$.

Therefore, we conclude that:
\[
\mathbf{V} \vee \mathbf{HSP}(W_2, Z_7, Z_8) = \widehat{\mathbf{V}}.
\]

\textbf{Case 5.} Suppose that $Z_2, W_2 \models u \approx v$ but $Z_7, Z_8 \not\models u \approx v$. Then the identity $u \approx v$ satisfies one of the following two conditions:
\begin{itemize}
    \item $m, n \geq 2$, there exists some $u_i \in u$ with $|u_i| = 1$ such that one of $\operatorname{occ}(u_i, u)$ and $\operatorname{occ}(u_i, v)$ is odd while the other is even. Moreover, there exists a nonempty subset $A \subseteq X$ such that the cardinalities
    \[
    \left| \left\{ i \mid u_i \in u,\ c(u_i) \subseteq A \right\} \right| \quad \text{and} \quad \left| \left\{ j \mid v_j \in v,\ c(v_j) \subseteq A \right\} \right|
    \]
    differ in parity;
    
    \item $m, n \geq 2$, there exists some $v_i \in v$ with $|v_i| = 1$ such that one of $\operatorname{occ}(v_i, v)$ and $\operatorname{occ}(v_i, u)$ is odd while the other is even. Moreover, there exists a nonempty subset $A \subseteq X$ such that the cardinalities
    \[
    \left| \left\{ i \mid u_i \in u,\ c(u_i) \subseteq A \right\} \right| \quad \text{and} \quad \left| \left\{ j \mid v_j \in v,\ c(v_j) \subseteq A \right\} \right|
    \]
    differ in parity.
\end{itemize}

In either case, the identity \( u \approx v \) implies an identity of the form
\[
x + u' \approx 2x + v',
\]
where \( x \notin u' \) and \( x \notin v' \).

By Lemma~\ref{3.64321}, it follows that
\[
\mathbf{V}_1 \subseteq \mathbf{V} \vee \mathbf{HSP}(Z_2, W_2).
\]

On the other hand, since \( \mathbf{V}_1 \models x + u' \approx 2x + v' \) but \( \mathbf{V}_1 \not\models x \approx 2x \), we conclude that \( Z_2, W_2 \in \mathbf{V}_1 \), and hence
\[
\mathbf{V} \vee \mathbf{HSP}(Z_2, W_2) \subseteq \mathbf{V}_1.
\]

Therefore, we obtain the equality:
\[
\mathbf{V}_1 = \mathbf{V} \vee \mathbf{HSP}(Z_2, W_2).
\]

\textbf{Case 6.} Suppose that $Z_2, Z_7 \models u \approx v$ but $W_2, Z_8 \not\models u \approx v$. Then the identity $u \approx v$ must satisfy one of the following three conditions:
\begin{itemize}
    \item $m = n = 1$, $|u_1|, |v_1| \geq 2$, and $c(u_1) \neq c(v_1)$;

    \item $m = 1$, $n \geq 2$, $|u_1| \geq 2$, and there exists a nonempty subset $A \subseteq X$ with $c(u_1) \subseteq A$ such that
    \[
    \left| \left\{ j \mid v_j \in v,\ c(v_j) \subseteq A \right\} \right| \text{ is even},
    \]
    and for any $v_j \in v$ with $|v_j| = 1$, we have $\operatorname{occ}(v_j, v)$ even;

    \item $n = 1$, $m \geq 2$, $|v_1| \geq 2$, and there exists a nonempty subset $A \subseteq X$ with $c(v_1) \subseteq A$ such that
    \[
    \left| \left\{ i \mid u_i \in u,\ c(u_i) \subseteq A \right\} \right| \text{ is even},
    \]
    and for any $u_j \in u$ with $|u_j| = 1$, we have $\operatorname{occ}(u_j, u)$ even.
\end{itemize}

In any of the above cases, the identity \( u \approx v \) implies the identity
\[
xy \approx 2xy.
\]

By Lemma~\ref{ccccc}, it follows that
\[
\mathbf{V}_1 \subseteq \mathbf{V} \vee \mathbf{HSP}(Z_2, Z_7).
\]

On the other hand, since \( \mathbf{V}_1 \models xy \approx 2xy \) and \( \mathbf{V}_1 \not\models x \approx 2x \), we conclude that \( Z_2, Z_7 \in \mathbf{V}_1 \), and hence
\[
\mathbf{V} \vee \mathbf{HSP}(Z_2, Z_7) \subseteq \mathbf{V}_1.
\]

Therefore, we conclude:
\[
\mathbf{V}_1 = \mathbf{V} \vee \mathbf{HSP}(Z_2, Z_7).
\]
\textbf{Case 7.} Suppose that $Z_2, Z_8 \models u \approx v$ but $W_2, Z_7 \not\models u \approx v$. Then the identity \( u \approx v \) must satisfy one of the following two conditions:
\begin{itemize}
    \item $m = 1$, $n \geq 2$, \( |u_1| \neq 1 \), and for all nonempty subsets \( A \subseteq X \) with \( c(u_1) \subseteq A \), we have
    \[
    \left| \left\{ i \mid v_i \in v,\ c(v_i) \subseteq A \right\} \right| \text{ is odd}, \quad
    \left| \left\{ i \mid v_i \in v,\ c(v_i) \not\subseteq A \right\} \right| \text{ is even},
    \]
    and there exists some \( v_j \in v \) such that \( |v_j| = 1 \) and \( \operatorname{occ}(v_j, v) \) is odd;

    \item $m \geq 2$, $n = 1$, \( |v_1| \neq 1 \), and for all nonempty subsets \( A \subseteq X \) with \( c(v_1) \subseteq A \), we have
    \[
    \left| \left\{ i \mid u_i \in u,\ c(u_i) \subseteq A \right\} \right| \text{ is odd}, \quad
    \left| \left\{ i \mid u_i \in u,\ c(u_i) \not\subseteq A \right\} \right| \text{ is even},
    \]
    and there exists some \( u_j \in u \) such that \( |u_j| = 1 \) and \( \operatorname{occ}(u_j, u) \) is odd.
\end{itemize}

\vspace{0.5em}

If \( N_2 \in \mathbf{V} \), then in the first (respectively, second) case, it must be that \( |v_j| \geq 2 \) for all \( v_j \in v \) (respectively, \( |u_i| \geq 2 \) for all \( u_i \in u \)), contradicting the existence of a block of length 1 with odd occurrence. Hence,
\[
\mathbf{V} \vee \mathbf{HSP}(Z_2, Z_8) = \mathbf{V} \vee \mathbf{HSP}(Z_2, W_2, Z_8)
\quad \text{or} \quad
\mathbf{V} \vee \mathbf{HSP}(Z_2, Z_8) = \mathbf{V} \vee \mathbf{HSP}(Z_2, Z_7, Z_8).
\]

\vspace{0.5em}

If \( N_2 \notin \mathbf{V} \), then the identity \( u \approx v \) implies the identity
\[
x^2 \approx x + 2x^2.
\]
By Lemma~\ref{3.35672}, it follows that
\[
\mathbf{V}_1 \subseteq \mathbf{V} \vee \mathbf{HSP}(Z_2, Z_8).
\]

On the other hand, since \( \mathbf{V}_1 \models x^2 \approx x + 2x^2 \) and \( \mathbf{V}_1 \not\models x \approx 2x \), we conclude that \( Z_2, Z_8 \in \mathbf{V}_1 \), and hence
\[
\mathbf{V} \vee \mathbf{HSP}(Z_2, Z_8) \subseteq \mathbf{V}_1.
\]

Therefore, we obtain:
\[
\mathbf{V}_1 = \mathbf{V} \vee \mathbf{HSP}(Z_2, Z_8).
\]

\textbf{Case 8.} Suppose that \( W_2, Z_7 \models u \approx v \) and \( Z_2, Z_8 \not\models u \approx v \). 

From \( Z_2, Z_8 \not\models u \approx v \), it follows that either $u = x$, and there exists a nonempty  $A \subseteq X$  with  $x \in A$  such that $\left| \{ j \mid v_j \in v, c(v_j) \subseteq A \} \right|$  is even,
or
\[
v = x, \quad \text{and } \left| \{ i \mid u_i \in u, c(u_i) \subseteq A \} \right| \text{ is even}.
\]

However, this contradicts the fact that \( W_2, Z_7 \models u \approx v \). 

Therefore, we conclude that
\[
\mathbf{V} \vee \mathbf{HSP}(W_2, Z_7) = \mathbf{V} \vee \mathbf{HSP}(Z_2, W_2, Z_7)
\quad \text{or} \quad
\mathbf{V} \vee \mathbf{HSP}(W_2, Z_7) = \mathbf{V} \vee \mathbf{HSP}(W_2, Z_7, Z_8).
\]

\textbf{Case 9.} Suppose that \( W_2, Z_8 \models u \approx v \) and \( Z_2, Z_7 \not\models u \approx v \). Then \( u \approx v \) satisfies one of the following two cases:
\begin{itemize}
  \item \( m = n = 1 \), \( |u_1| = 1 \), \( |v_1| \geq 2 \), and \( c(v_1) = c(u_1) \);
  \item \( m = n = 1 \), \( |v_1| = 1 \), \( |u_1| \geq 2 \), and \( c(u_1) = c(v_1) \).
\end{itemize}

If \( N_2, T_2 \notin \mathbf{V} \), then in either case \( u \approx v \) can imply the identity \( x \approx x^2 \). By Lemma \ref{3.4}, \(\mathbf{V}_1\) is a subvariety of \(\mathbf{V} \vee \mathbf{HSP}(W_2, Z_8)\). On the other hand, since \(\mathbf{V}_1 \models x \approx x^2\) and \(\mathbf{V}_1 \not\models x \approx x + x\), it follows that \( W_2, Z_8 \in \mathbf{V}_1 \), and hence \(\mathbf{V} \vee \mathbf{HSP}(W_2, Z_8)\) is a subvariety of \(\mathbf{V}_1\). Thus,
\[
\mathbf{V}_1 = \mathbf{V} \vee \mathbf{HSP}(W_2, Z_8).
\]

If \( N_2 \in \mathbf{V} \), then by Lemma \ref{1.1} (5), we must have \( |u_1| = |v_1| = 1 \), which is a contradiction. Therefore,
\[
\mathbf{V} \vee \mathbf{HSP}(W_2, Z_8) = \mathbf{V} \vee \mathbf{HSP}(Z_2, W_2, Z_8) \quad \text{or} \quad \mathbf{V} \vee \mathbf{HSP}(W_2, Z_8) = \mathbf{V} \vee \mathbf{HSP}(W_2, Z_7, Z_8).
\]

If \( T_2 \in \mathbf{V} \), then by Lemma \ref{1.1} (6), we must have \( |u_1|, |v_1| \geq 2 \), which contradicts the conditions above. Hence,
\[
\mathbf{V} \vee \mathbf{HSP}(W_2, Z_8) = \mathbf{V} \vee \mathbf{HSP}(Z_2, W_2, Z_8) \quad \text{or} \quad \mathbf{V} \vee \mathbf{HSP}(W_2, Z_8) = \mathbf{V} \vee \mathbf{HSP}(W_2, Z_7, Z_8).
\]

\textbf{Case 10.} Suppose \( Z_7, Z_8 \models u \approx v \) and \( Z_2, W_2 \not\models u \approx v \). Then \( u \approx v \) satisfies one of the following two cases:
\begin{itemize}
  \item \( m = 1, n \geq 2 \), \( |u_1| = 1 \) and \( occ(u_1, v) \) is odd, and for all nonempty subsets \( A \subseteq X \) with \( c(u_1) \subseteq A \), the number
  \[
  |\{ j \mid v_j \in v,\, c(v_j) \subseteq A \}|
  \]
  is odd;
  
  \item \( n = 1, m \geq 2 \), \( |v_1| = 1 \) and \( occ(v_1, u) \) is odd, and for all nonempty subsets \( A \subseteq X \) with \( c(v_1) \subseteq A \), the number
  \[
  |\{ i \mid u_i \in u,\, c(u_i) \subseteq A \}|
  \]
  is odd.
\end{itemize}

If \( T_2 \notin \mathbf{V} \), then in either case \( u \approx v \) can imply the identity
\[
x \approx x + 2x^2.
\]
By Lemma \ref{3.555}, \(\mathbf{V}_1\) is a subvariety of \(\mathbf{V} \vee \mathbf{HSP}(Z_7, Z_8)\). On the other hand, since \(\mathbf{V}_1 \models x \approx x + 2x^2\) and \(\mathbf{V}_1 \not\models x \approx 2x\), it follows that \(Z_7, Z_8 \in \mathbf{V}_1\), and thus \(\mathbf{V} \vee \mathbf{HSP}(Z_7, Z_8)\) is a subvariety of \(\mathbf{V}_1\). Therefore,
\[
\mathbf{V}_1 = \mathbf{V} \vee \mathbf{HSP}(Z_7, Z_8).
\]

If \( T_2 \in \mathbf{V} \), then by Lemma \ref{1.1} (6), there exist blocks \( u_i \in u \) and \( v_i \in v \) such that \( |u_i| \geq 2 \) and \( |v_i| \geq 2 \), contradicting the conditions above. Hence,
\[
\mathbf{V} \vee \mathbf{HSP}(Z_7, Z_8) = \mathbf{V} \vee \mathbf{HSP}(Z_2, Z_7, Z_8) \quad \text{or} \quad \mathbf{V} \vee \mathbf{HSP}(Z_7, Z_8) = \mathbf{V} \vee \mathbf{HSP}(W_2, Z_7, Z_8).
\]

\textbf{Case 11.} Suppose \( Z_2 \models u \approx v \) and \( W_2, Z_7, Z_8 \not\models u \approx v \). Then \( u \approx v \) satisfies one of the following two cases:
\begin{itemize}
  \item \( m = 1, n \geq 2 \), \( |u_1| \neq 1 \), and there exists a nonempty subset \( A \subseteq X \) with \( c(u_1) \subseteq A \) such that
  \[
  |\{ i \mid v_i \in v,\, c(v_i) \subseteq A \}|
  \]
  is even, and there exists \( v_j \) with \( |v_j| = 1 \) such that \( occ(v_j, v) \) is odd;
  
  \item \( m \geq 2, n = 1 \), \( |v_1| \neq 1 \), and there exists a nonempty subset \( A \subseteq X \) with \( c(v_1) \subseteq A \) such that
  \[
  |\{ i \mid u_i \in u,\, c(u_i) \subseteq A \}|
  \]
  is even, and there exists \( u_j \) with \( |u_j| = 1 \) such that \( occ(u_j, u) \) is odd.
\end{itemize}

If \( N_2 \in \mathbf{V} \), then for the first (second) case, we have \( |v_j| \geq 2 \) (\( |u_i| \geq 2 \)), contradicting the assumption. Thus,
\[
\mathbf{V} \vee \mathbf{HSP}(Z_2) = \mathbf{V} \vee \mathbf{HSP}(Z_2, Z_7)\] or \[\mathbf{V} \vee \mathbf{HSP}(Z_2) = \mathbf{V} \vee \mathbf{HSP}(Z_2, W_2)\] or \[\mathbf{V} \vee \mathbf{HSP}(Z_2) = \mathbf{V} \vee \mathbf{HSP}(Z_2, Z_8).
\]

If \( N_2 \notin \mathbf{V} \), then \( u \approx v \) can imply the identity
\[
x^2 \approx x + x^2.
\]
By Lemma \ref{3.3567}, \(\mathbf{V}_1\) is a subvariety of \(\mathbf{V} \vee \mathbf{HSP}(Z_2)\). On the other hand, since \(\mathbf{V}_1 \models x^2 \approx x + x^2\) and \(\mathbf{V}_1 \not\models x \approx 2x\), it follows that \( Z_2 \in \mathbf{V}_1 \), and hence
\[
\mathbf{V} \vee \mathbf{HSP}(Z_2) \subseteq \mathbf{V}_1.
\]
Therefore,
\[
\mathbf{V}_1 = \mathbf{V} \vee \mathbf{HSP}(Z_2).
\]

\textbf{Case 12.} Suppose \( W_2 \models u \approx v \) and \( Z_2, Z_7, Z_8 \not\models u \approx v \).

From \( W_2 \models u \approx v \) and \( Z_2 \not\models u \approx v \), it follows that \( u \approx v \) satisfies one of the following two cases:
\begin{itemize}
  \item \( m = n = 1 \), \( |u_1| = 1 \), \( |v_1| \geq 2 \), and \( c(v_1) = c(u_1) \);
  \item \( m = n = 1 \), \( |v_1| = 1 \), \( |u_1| \geq 2 \), and \( c(u_1) = c(v_1) \).
\end{itemize}

It follows that \( Z_8 \models u_1 \approx v_1 \), a contradiction.

Thus,
\[
\mathbf{V} \vee \mathbf{HSP}(W_2) = \mathbf{V} \vee \mathbf{HSP}(Z_2, W_2)
\]
or
\[
\mathbf{V} \vee \mathbf{HSP}(W_2) = \mathbf{V} \vee \mathbf{HSP}(W_2, Z_7)
\]
or
\[
\mathbf{V} \vee \mathbf{HSP}(W_2) = \mathbf{V} \vee \mathbf{HSP}(W_2, Z_8).
\]

\textbf{Case 13.} Suppose \( Z_7 \models u \approx v \) and \( Z_2, W_2, Z_8 \not\models u \approx v \). Then \( u \approx v \) satisfies one of the following two cases:
\begin{itemize}
  \item \( m = 1, n \geq 2 \), \( |u_1| = 1 \), 
  \[
    \{ v_j \mid |v_j| = 1, \, \mathrm{occ}(v_j, v) \text{ is odd} \} = \{ u_1 \},
  \]
  and there exists a nonempty subset \( A \subseteq X \) with \( c(u_1) \subseteq A \) such that 
  \[
    |\{ j \mid v_j \in v, c(v_j) \subseteq A \}|
  \]
  is even;
  \item \( n = 1, m \geq 2 \), \( |v_1| = 1 \), 
  \[
    \{ u_i \mid |u_i| = 1, \, \mathrm{occ}(u_i, u) \text{ is odd} \} = \{ v_1 \},
  \]
  and there exists a nonempty subset \( A \subseteq X \) with \( c(v_1) \subseteq A \) such that
  \[
    |\{ i \mid u_i \in u, c(u_i) \subseteq A \}|
  \]
  is even.
\end{itemize}

If \( T_2 \notin \mathbf{V} \), then in either case \( u \approx v \) implies the identity \( x \approx x + x^2 \). By Lemma \ref{3.55}, \(\mathbf{V}_1\) is a subvariety of \(\mathbf{V} \vee \mathbf{HSP}(Z_7)\). On the other hand, since \(\mathbf{V}_1 \models x \approx x + x^2\) and \(\mathbf{V}_1 \not\models x \approx 2x\), it follows that \( Z_7 \) is a member of \(\mathbf{V}_1\). Hence, \(\mathbf{V} \vee \mathbf{HSP}(Z_7)\) is a subvariety of \(\mathbf{V}_1\). Thus,
\[
\mathbf{V}_1 = \mathbf{V} \vee \mathbf{HSP}(Z_7).
\]

If \( T_2 \in \mathbf{V} \), then by Lemma \ref{1.1} (6), we have
\[
\{ u_i \in u \mid |u_i| \geq 2 \} \neq \emptyset, \quad \text{and} \quad \{ v_i \in v \mid |v_i| \geq 2 \} \neq \emptyset,
\]
which is a contradiction.

Thus,
\[
\mathbf{V} \vee \mathbf{HSP}(Z_7) = \mathbf{V} \vee \mathbf{HSP}(Z_2, Z_7)
\]
or
\[
\mathbf{V} \vee \mathbf{HSP}(Z_7) = \mathbf{V} \vee \mathbf{HSP}(W_2, Z_7)
\]
or
\[
\mathbf{V} \vee \mathbf{HSP}(Z_7) = \mathbf{V} \vee \mathbf{HSP}(Z_7, Z_8).
\]

\textbf{Case 14.} Suppose \( Z_8 \models u \approx v \) and \( Z_2, W_2, Z_7 \not\models u \approx v \). Then \( u \approx v \) satisfies one of the following two cases:
\begin{itemize}
  \item \( m = 1, n \geq 2 \), \( |u_1| = 1 \), 
  \[
    (\forall \emptyset \neq A \subseteq X, c(u_1) \subseteq A) \quad |\{ j \mid v_j \in v, c(v_j) \subseteq A \}| \text{ is odd},
  \]
  and 
  \[
    \{ v_j \mid |v_j| = 1, \text{occ}(v_j, v) \text{ is odd} \} \neq \{ u_1 \};
  \]
  \item \( n = 1, m \geq 2 \), \( |v_1| = 1 \), 
  \[
    (\forall \emptyset \neq A \subseteq X, c(v_1) \subseteq A) \quad |\{ i \mid u_i \in u, c(u_i) \subseteq A \}| \text{ is odd},
  \]
  and 
  \[
    \{ u_i \mid |u_i| = 1, \text{occ}(u_i, u) \text{ is odd} \} \neq \{ v_1 \}.
  \]
\end{itemize}

If \( T_2 \notin \mathbf{V} \), then in either case \( u \approx v \) implies the identity
\[
x \approx 2x + x^2.
\]
By Lemma \ref{3.5555}, \(\mathbf{V}_1\) is a subvariety of \(\mathbf{V} \vee \mathbf{HSP}(Z_8)\). On the other hand, since
\[
\mathbf{V}_1 \models x \approx 2x + x^2 \quad \text{and} \quad \mathbf{V}_1 \not\models x \approx 2x,
\]
it follows that \( Z_8 \in \mathbf{V}_1 \), so
\[
\mathbf{V} \vee \mathbf{HSP}(Z_8) \subseteq \mathbf{V}_1.
\]
Hence,
\[
\mathbf{V}_1 = \mathbf{V} \vee \mathbf{HSP}(Z_8).
\]

If \( T_2 \in \mathbf{V} \), then we have 
\[
(\forall v_j) \quad |c(v_j)| = 1 \quad (\text{resp. } (\forall u_i) \quad |c(u_i)| = 1).
\]
Since \( Z_8 \models u \approx v \), it follows that
\[
(\forall A \subseteq X, c(u_1) \subseteq A) \quad |\{ j \mid v_j \in v, c(v_j) \subseteq A \}| \]
\[\quad \text{(resp. }(\forall A \subseteq X, c(v_1) \subseteq A) \quad |\{ i \mid u_i \in u, c(u_i) \subseteq A \}| \text{)} \quad \text{is odd}.
\]
This implies
\[
Z_7 \models u \approx v,
\]
a contradiction.

Thus,
\[
\mathbf{V} \vee \mathbf{HSP}(Z_8) = \mathbf{V} \vee \mathbf{HSP}(Z_2, Z_8)
\]
 or
\[
\mathbf{V} \vee \mathbf{HSP}(Z_8) = \mathbf{V} \vee \mathbf{HSP}(W_2, Z_8)
\]
or
\[
\mathbf{V} \vee \mathbf{HSP}(Z_8) = \mathbf{V} \vee \mathbf{HSP}(Z_7, Z_8).
\]

Based on the 14 cases discussed above, we analyze the four possible configurations of the elements \( N_2 \) and \( T_2 \) with respect to their membership in the set \( \mathbf{V} \). These configurations are classified as follows:

\begin{equation*}
\begin{aligned}
&\text{Case 1:} \quad N_2, T_2 \in \mathbf{V} && \text{(both \( N_2 \) and \( T_2 \) are members of \( \mathbf{V} \))}. \\
&\text{Case 2:} \quad N_2 \notin \mathbf{V},\ T_2 \in \mathbf{V} && \text{(only \( T_2 \) is a member of \( \mathbf{V} \))}. \\
&\text{Case 3:} \quad N_2 \in \mathbf{V},\ T_2 \notin \mathbf{V} && \text{(only \( N_2 \) is a member of \( \mathbf{V} \))}.\\
&\text{Case 4:} \quad N_2, T_2 \notin \mathbf{V} && \text{(neither \( N_2 \) nor \( T_2 \) is a member of \( \mathbf{V} \))}.
\end{aligned}
\end{equation*}

These four membership patterns correspond to distinct configurations within the lattice structure defined over the interval \( [\mathbf{V}, \widehat{\mathbf{V}}] \), as illustrated in Figure 1. Each configuration reflects a unique arrangement of elements under the partial ordering induced by the lattice.

\end{Proof}

\begin{Theorem} \label{3.6}
    ${\cal L}(\mathbf{Sr}(2))$ is a lattice of order \emph{480}.
\end{Theorem}

\begin{Proof}
    By Equation~(\ref{24}) and Theorem \ref{3.5}, we conclude that ${\cal L}(\mathbf{Sr}(2))$ contains exactly $480$ $(2^4*(6+8+9+7)=480)$ elements.
\end{Proof}
\begin{Lemma}\cite{rz} \label{3.7}
    If a variety $\mathbf{V}$ is finitely based and ${\cal L}(\mathbf{V})$ is a finite lattice, then $\mathbf{V}$ is hereditarily finitely based 
\end{Lemma}

By Theorem \ref{3131}, \ref{3.6} and Lemma \ref{3.7}, we now immediately deduce the following:

\begin{Theorem}
    $\mathbf{Sr}(2)$ is hereditarily finitely based.
\end{Theorem}

\vskip 15pt
\textbf{Acknowledgements}
	The research is supported  by  by National Natural Science Foundation of China (12371024) and  Science and Technology Research Program of Chongqing Education Commission of China (KJZD-K202401102).
\vskip 15pt
\begin{flushleft}
\textbf{Conflicts of interest} All authors declare that they have no conflicts of interest.
\end{flushleft}


\begin{thebibliography}{aa}
\bibitem{bur} Burris, S., Sankappanavar, H.P., A Course in Universal Algebra,
Springer, New York (1981).
 \bibitem{dlk} Dolinka I., A nonfintely based finite semiring, Int. J. Algebra Comput., 2007, 17(8): 1537-1551.
\bibitem{bk} El Bashir, R., Kepka, T., Congruence-simple semirings, Semigroup Forum, 2007, 75: 588-608.
\bibitem{gho} Ghosh, S., Pastijn, F., Zhao, X.Z., Varieties Generated by Ordered Bands I, Order, 2005, 22: 109-128.
\bibitem{glk} G{\l}azek, K., A guide to the literature on semirings and their applications in mathematics and information science,
Kluwer Academic Publishers, Dordrecht (2001).
\bibitem{gol} Golan, J.S., The theory of semirings with applications in mathematics and theoretical computer science,
Longman Scientific and Thechnical, Harlow (1992).
\bibitem{gf} Guzm$\acute{a}$n, F., The variety of Boolean Semirings,
J. Pure Appl. Algebra, 1992, 8: 253-270.
\bibitem{hw}  Hebisch, U., Weinert, H.J., Semirings. Algebraic Theory and Applications in Computer Science,
World Scientific, Singapore (1998).
\bibitem{jmh} Howie, J. M., Fundaments of Semigroup Theory, Oxford Science Publication, Oxford, 1995.
Order, 2005, 22: 109-128.
\bibitem{jrz}Jackson, M., Ren, M.M., Zhao, X.Z., Nonfinitely based ai-semirings with finitely based semigroup reducts, Journal of Algebra, 2022, 611: 211-245.
 \bibitem{krrl} Kruse R. L., Identities satisfied by a finite ring, J. Algebra, 1973, 26(2): 298-318.
\bibitem{lci} Lyndon R. C., Identities in two-valued calculi, Trans. Amer. Math. Soc., 1951, 71(3): 457-457.
\bibitem{lcii} Lyndon R. C., Identities in finite algebras, Proc. Amer. Math. Soc., 1954, 5: 8-9.
 \bibitem{liv}  L'vov I. V., Varieties of associative rings, I. Algebra and Logic, 1973, 12(3): 150-167.
 \bibitem{mcnr} McKenzie R., Equational bases for lattice theories, Math. Scand., 1970, 27: 24-38.
\bibitem{mker}McKenzie R., Tarski's finite basis problem is undecidable, Int. J. Algebra Comput., 1996, 6(1): 49-104.
\bibitem{mgfw}McNulty G. F., Willard R., The Chautauqua Problem, Tarski's Finite Basis Problem, and Residual
Bounds for 3-element Algebras. In progress.
\bibitem{oats} Oates S., Powell M. B., Identical relations in finite groups, J. Algebra, 1964, 1(1): 11-39.
\bibitem{perp} Perkins P., Bases for equational theories of semigroups, J. Algebra, 1969, 11(2): 298-314.
\bibitem{pet} Petrich, M., Reilly, N.R., Completely Regular Semigroups, Wiley, New York, (1999).
\bibitem{259} Polin, S.V., Minimal varieties of semirings, Matematicheskie Zametki, 1980, 27(4): 527-537.
\bibitem{rz}  Ren, M.M., Zeng, L.L., On a hereditarily finitely based ai-semiring variety, Soft Computing, 2019, 23(16): 6819-6825.
\bibitem{srii} Shao, Y., Ren, M.M., On the variety generated by all ai-semirings of order two,
Semigroup Forum, 2015, 91: 171-184.
\bibitem{tasa}Tarski A., Equational logic and equational theories of algebras, Stud. Logic Found. Math., 1968, 50: 275-288.
\bibitem{vp} Vechtomov, E.M., Petrov, A.A., Multiplicatively idempotent semirings, Journal of Mathematical Sciences, 2015, 206(6): 634-653.
\bibitem{rzcsd} Zhao, X.Z., Ren, M.M., Crvenkovi\'{c}, S., Shao, Y., Dapi\'{c}, P., The variety generated by semirings of order three, Ural Mathematical Journal, 2020, 6(2): 117-132.
\end{thebibliography}
\end{document}